% plain TEX
\magnification=\magstep1
\dimen100=\hsize
\parskip= 6pt

\font\ninerm=cmr10 at 10truept
\font\eightrm=cmr8
\font\svntnrm = cmr17
\font\sc=cmcsc10
\font \eightsc = cmcsc8

\font\tenmsy=msbm10
\font\sevenmsy=msbm7
\font\fivemsy=msbm5
\newfam\msyfam
\def\msy{\fam\msyfam\tenmsy}
\textfont\msyfam=\tenmsy
\scriptfont\msyfam=\sevenmsy
\scriptscriptfont\msyfam=\fivemsy

\def\bbc{{\msy C}}

\def\bbp{{\msy P}}

\def\bbr{{\msy R}}

\def\bbz{{\msy Z}}

\def\bk{\bigskip}

\def\ni{\noindent}

\font\aa=eufm10

\def\bk{\bigskip}
\def\Got#1{\hbox{\aa#1}}

\def\Ra*{(R_{a})_{*}}
\def\ada-1{ad_{a^{-1}}~}
\def\p-1u{\pi^{-1}(U)}

\def\liet{{ \Got t}}
\def\lieg{{ \Got g}}
\def\mt {{\cal M}^\liet}
\def\mg {{\cal M}^\lieg}
\def\mT {{\cal M}^T}
\def\mG {{\cal M}^G}
\def\alc{{\Delta}}
\def \Epsilon{ {\cal E} } 
\input amssym.def
\input amssym.tex

\centerline{\svntnrm  Representations with Weighted Frames} 

\centerline{\svntnrm  and Framed Parabolic Bundles}
\medskip

\bigskip
\centerline{\sc   J.C. Hurtubise~~ L.C. Jeffrey\footnote{}{\ninerm
During the preparation of this work the   authors were
supported    NSERC and FCAR grants.}}

\bigskip
\centerline{\vbox{\hsize = 5.85truein
\baselineskip = 12.5truept
\eightrm {\eightsc Abstract:} There is a well-known correspondence between the symplectic 
variety $M_h$ of
representations of
 the fundamental group of a punctured Riemann surface into a compact 
Lie group $G$, 
with fixed conjugacy
classes $h$ at the punctures, and a complex variety ${\cal M}_h$ of 
holomorphic bundles on the unpunctured surface with 
a parabolic structure at the puncture points. For $G = SU(2)$, we build a symplectic variety $P$ of pairs 
(representations of the fundamental group into $G$, ``weighted frame'' at the puncture points),
and a corresponding complex variety ${\cal P}$ of moduli of ``framed parabolic bundles'', 
which encompass 
respectively all of the spaces $M_h$, ${\cal M}_h$, in the sense that one can obtain 
$M_h$ from $P$ by symplectic reduction, and ${\cal M}_h$ from ${\cal P}$ 
by a complex quotient. This allows us to explain certain features of the toric geometry 
of the $SU(2)$ moduli spaces discussed by Jeffrey and Weitsman, by giving the actual toric 
variety associated with their integrable system.}}

{\bf 1. Introduction}

The geometry of moduli spaces of vector bundles over a Riemann surface
exhibits in a beautiful way the interplay between two types of structure.
On one hand, one can think of these spaces as classifying holomorphic
structures, and the emphasis is on complex geometry. On the other, via the
Narasimhan-Seshadri theorem [NS], one can also represent the moduli space as 
a  space of flat unitary connections on the surface, or, integrating, in terms of 
representations of the fundamental group into unitary groups; here the emphasis is
on the symplectic geometry. The key link between the two is 
 provided by the fact that on both sides one can view the moduli space
as being obtained by a quotient construction, starting from the same 
infinite dimensional space of unitary connections [AB]. On the holomorphic side, one 
quotients by the action of the group of complex gauge transformations; on the 
symplectic side, one performs a symplectic reduction, first setting the curvature (moment
map) to zero, then quotienting out by the group of unitary gauge
 transformations.
The general equivalence between geometric invariant theory quotients 
and symplectic quotients tells us that we should obtain the same result.

Following Mehta and Seshadri,
one can extend this picture to the case of a Riemann surface with punctures; 
on the symplectic side, one considers a moduli space of flat connections
 over the surface, such that the conjugacy class of the holonomy of
 the connection
around the punctures is fixed;  on the holomorphic side, the corresponding
object is a moduli space of vector bundles over the filled-in Riemann surface, with the 
extra data of a parabolic structure at the marked points, that is,  a
flag of subspaces in the fiber of the bundle at each of the marked points. The 
conjugacy class determines the choice of polarisation that one uses to define
 the quotient.

 From the symplectic point of view, it is tempting to try and obtain the moduli spaces
of flat connections on a surface by ``glueing'' moduli spaces of connections
over the elements of a decomposition of the surface into punctured surfaces. 
This requires moduli spaces over the punctured surface which are larger than the
parabolic moduli; one needs, at the very least, to have all possible conjugacy
 classes of the holonomy around the puncture, as well as some sort of framing at the
punctures to remove the ambiguity in glueing. The smallest moduli spaces with  which one could
have hoped to
achieve this  are the $T$-extended moduli spaces of Jeffrey [J1], where $T$ is 
the maximal torus for our group $G$. 
These are spaces 
of equivalence classes of 
flat connections on the surface, along with a trivialisation at each puncture in 
which the holonomy at the puncture lies in the maximal torus $T$. The holonomy
 variables (eigenvalues of the holonomy)
are the ``associated momenta" of the framing variables, so that 
if one decomposes a surface $\Sigma$ into two punctured surfaces $\Sigma_1$,
$\Sigma_2$, the glueing operation on the moduli spaces of
matching the holonomies and then quotienting out the equivalent framings amounts to 
a symplectic quotient of the product of the two  moduli spaces associated to 
 $\Sigma_1$,$\Sigma_2$ by the diagonal $T$-action on the framings. This glueing
works well
at the points for which the holonomy is regular, but fails when it is not; for one, the extended
moduli are no longer necessarily symplectic at these points [J1]; also, the quotient space
of the moduli associated to the two surfaces maps to the moduli of the glued surface, but 
is no longer an isomorphism over the connections with non-regular holonomy along the 
circle of glueing. 
Various schemes can be employed to repair  this problem, for example the $q$-Hamiltonian spaces of 
Alekseev, Malkin and Meinrenken [AMM], or the loop group valued framings of Meinrenken 
and Woodward
[MW].

Our aim is different. For $G= SU(2), T= S^1$, we  obtain a moduli space $P$
which 
is symplectic, and which encompasses all of the parabolic moduli, in some sense, 
``repairing'' the $T$-extended moduli space. This new space will have a Hamiltonian
$T$-action associated to each puncture, whose moment map is the holonomy, and whose symplectic $T$-quotients
will be the parabolic moduli spaces. The moduli space gets interpreted as a moduli 
space of bundles with ``weighted frames'' at the punctures; see section 2. In section 3 we will
see that the spaces associated to two punctured surfaces can indeed be glued, though 
one does not obtain the moduli space of bundles over the glueing of the two surfaces. 

The space $P$  allows one to clear up a question concerning the Goldman flows [G1]
 on the moduli space. When one is dealing with the product of the two $T$-
extended moduli spaces associated to two punctured surfaces $\Sigma_1$ and $\Sigma_2$, one has an action
of $T\times T$, given by the action of $T$ on the trivialisations over 
the two punctures. Let $\Sigma$ be the surface obtained by glueing $\Sigma_1$ 
and $\Sigma_2$, identifying curves around the punctures to obtain a single curve 
 $C$ in $\Sigma$. In the corresponding glueing of the $T$-extended moduli 
spaces,  we are just quotienting
out the anti-diagonal $T$, and so the group $1\times T$ acts on the quotient. This
action only descends
to the moduli of the glued surface over the set of connections with regular holonomy
along $C$ and gives a 
densely defined $T$-flow on the moduli space of the glued surface $\Sigma$, which 
is ill-defined precisely when the holonomy along $C$ is not regular.
In the case of $G=SU(2), T= S^1$, doing this for the $3g-3$ boundary 
circles of a decomposition of the genus $g$ surface $\Sigma$, gives a Lagrangian
foliation of the moduli space corresponding to a densely defined $(S^1)^{3g-3}$
action, with moment map onto a polytope $\Gamma$, and ``almost'' 
gives the moduli space the structure of a toric variety. This 
``almost'' is tantalizing: if the variety were toric, one could
get a proof of the Verlinde formulae for the dimensions of the 
spaces of sections of line bundles over the moduli space, simply by counting lattice points in 
the moment polytope; indeed, 
 Jeffrey and Weitsman, in [JW], count these points,
   and  obtain the answers 
predicted by the Verlinde formulae. If the moduli space had only been toric, this would have been
a proof, and indeed, understanding how the structure fails to be toric was one of our main
motivations for this paper. Our spaces $P$ provide the answer to this question.
One can glue the spaces $P_1, P_2 $ correponding to $\Sigma_1, \Sigma_2$
by taking the anti-diagonal $T$-quotient: 
 one obtains a moduli space $\hat P$ associated to the glueing
of $\Sigma_1$ and $\Sigma_2$. If one does this for a
 family of trinions associated to a decomposition $D$
of a Riemann surface,
one obtains a toric variety $P^D$ 
whose moment polytope is precisely the 
$\Gamma$ of the Goldman flows; the moduli space $\mT$ maps to $P^D$
in a generically bijective fashion, and the Goldman flows are the lifts of the flows on the 
toric variety.

One could also  hope that the space $P$ corresponds, as in the Narasimhan-Seshadri theorem, to a
moduli space of holomorphic objects over the curve $\Sigma$, where one no longer has punctures but 
simply marked points. We will show that this is also the case: the objects 
one considers are framed parabolic bundles, parabolic bundles with the extra structure
of a volume form  on each of the successive quotients of the vector spaces in the flag.
This space is not complete, and to complete it, one must allow the bundles to degenerate into certain 
sheaves with a small amount of torsion. Our construction, which we carry out
 here for the 
$Sl(2, \bbc)$-case only, involves using the ``master space'' construction of Thaddeus in an
essential way. One can again glue in the holomorphic category the moduli spaces $P_1$, 
$P_2$ associated to $\Sigma_1, \Sigma_2$ by what is now an anti-diagonal $T_\bbc$ geometric invariant theory
quotient, and the holomorphic objects one obtains live over the singular curve obtained 
by identifying the marked points of $\Sigma_1, \Sigma_2$. 

In section 2, we construct our moduli space $P$, from a symplectic point of view.
Section three is devoted to the glueing construction and
 the example of  the moduli space over the three-punctured sphere.
In section 4, we build the holomorphic versions of our spaces.

We would like to thank Reyer Sjamaar for some crucial conversations.

\noindent{\bf 2. Extended moduli spaces, and implosions}

Let $G$ be a  simply connected compact Lie group, with maximal torus $T$,
and let $\lieg, \liet$ be their Lie algebras. We define  
the {\it{open fundamental alcove}} as the space 
$$ \Delta^0 = \{ \xi 
\in \liet : (\alpha_j, \xi)>0, ~\alpha_{\rm max} (\xi) < 1 \}\eqno (2.1)$$
where $\alpha_{\rm max}$ is the highest root and the $\alpha_j$ 
are the simple roots.
The { \it fundamental alcove} is the closure of $\Delta^0$ in 
$\liet$  and will be  denoted $\Delta$. Every point in $G$ is conjugate to a unique
element in $\exp \Delta$. The regular elements of $G$ correspond to the 
interior of $\Delta$. For $G= SU(2)$, the fundamental alcove is an interval,
representing the elements $diag(e^{i\pi\gamma}, e^{-i\pi\gamma}), \ \gamma\in 
[0,1]$ in $SU(2)$.

In what follows, let $\Sigma$ denote a Riemann surface of genus $g$
 with punctures at points $p_1,..,p_n, n\ge 0$. 

We want to describe the surface as a polygon with sides identified. 
We can choose closed paths $a_1, ..,a_g,b_1,..,b_g, c_1,..,c_n$ representing 
 cycles of a homology basis of 
$\Sigma$. The $a_i$, $b_i$ originate at a point $\tilde p_1$ near the puncture 
$p_1$, and represent the standard homology basis of the closure of the surface. 
The $c_i$ originate at points $\tilde p_i$ near $p_i$, and are small loops surrounding the punctures.
 As well, 
 we can pick paths $d_i, i=2,..,n$, joining $p_1$ to $p_i$.
 We can further assume that the paths do not 
intersect except at their end points, and are such  that we can write $\Sigma $ in the usual way as
a polygon with sides (in order) 
$$a_1,b_1,a_1^{-1},b_1^ {-1}, a_2, b_2,...,a_g^{-1},b_g^{-1},c_1, d_2,c_2 ,
d_2^{-1},d_3,c_3, d_3^{-1}, ...,c_n,d_n^{-1}\eqno (2.2)$$
with  identifications, as given, along the edges. 
The fundamental group of $\Sigma$
is then the free group on the generators $a_i, b_i, c_1, d_jc_jd_j^{-1}$ with the relation 
that the product $a_1b_1a_1^{-1}\cdots d_n^{-1}$ along the sides of the polygon is one.

We next define some moduli spaces. If $(A_1, \dots, A_g, B_1,.., B_g, ,D_2,.., 
D_n, ,C_1,.., C_n)\in G^{2g+ 2n-1} $, define
$$ \eqalign{ F: G^{2g+ 2n-1} &\rightarrow G\cr 
 (A_i, B_i, D_j, C_j)& \mapsto 
A_1B_1A_1^{-1}B_1^{-1}\cdots A_g^{-1}B_g^{-1}C_1D_2C_2D_2^{-1}\cdots C_nD_n^{-1}. }
\eqno (2.3)$$%{(*mapdef)}
 The {\it $\liet$-extended moduli space} 
$\mt$ is defined by 
$$\mt =\{ (A_i,B_i,D_j,\gamma_j)\in { G}^{2g+n-1}
\times \liet^n :F(A_i, B_i, D_j, C_j) =1 \}\eqno (2.4)$$
setting  $C_i=exp(\gamma_i)$.
There is a corresponding {\it $\lieg$-extended moduli space}
 where the $\gamma_j$ are not restricted
to lie in $\liet$ but can lie anywhere in $\lieg$:
$$\mg =\{ (A_i,B_i,D_j,\gamma_j)\in { G}^{2g+n-1}
\times \lieg^n : 
F(A_i, B_i, D_j, C_j) =1\}\eqno (2.5)$$
where as before $C_i = \exp (\gamma_i)$.
One may also define the {\it $T, G$- extended moduli spaces} 
$$\mT =\{ (A_i,B_i,D_j,C_j )\in { G}^{2g+n-1}
\times (\exp \alc)^n : 
F(A_i, B_i, D_j, C_j) =1 \}\eqno (2.6)$$
and
$$\mG =\{ (A_i,B_i,D_j,C_j )\in { G}^{2g+n-1}
\times G^n : 
F(A_i, B_i, D_j, C_j) =1\}\eqno (2.7)$$
There is an  obvious map from $\mg$ to $\mG$ (sending
the point $(A_j, B_j, D_j, \gamma_j)$ to the point
$(A_j, B_j, D_j, C_j)$).

Finally, one may make a minor modification of  these spaces: choosing an 
element $c \in Z(G)$ we may define 
$$\mg_c 
 =\{ (A_i,B_i,D_j,\gamma_j)\in { G}^{2g+n-1}
\times \lieg^n : 
A_1B_1A_1^{-1}\cdots B_g^{-1}C_1D_2C_2D_2^{-1}\cdots C_nD_n^{-1}=c \}\eqno (2.8)$$
and $\mt_c$, $\mG_c$ and $\mT_c$ are defined in the same way, replacing $1$ by $c$
in the definitions.

 These spaces can be thought of as moduli 
spaces of flat connections on $\Sigma$, with some additional  structure at the 
marked points. The space $\mt$ is the quotient of a space of flat connections on a Riemann 
surface $\Sigma$ with $n$ boundary components, where we have imposed some boundary conditions
near the punctures, in particular
fixing the connection form around each boundary component to lie in 
$\liet$. We divide only by those gauge transformations which are 
trivial on the boundary. The space $\mg$ is defined similarly, removing the  
constraint that the form lie in $\liet$ near the boundary. (See [J1].)
The spaces $\mT, \mG$  can similarly be obtained 
from a space of flat connections by quotienting by a suitable space of 
gauge transformations which are trivial only at marked points near the boundaries.
In other words, at a point near each
of the boundary components we have fixed a framing or trivialization
of the principal bundle on which the connection is defined. The space $\mT$
is a quotient of the space $\mt$, by the $n$-fold product $(W\ltimes\bbz^r)^n$
of the semidirect product of the Weyl group $W$  and the $\bbz^r$ of translations
such that $T=\liet/\bbz^r$. Alternately, $\mT$ can be thought of as a subspace
of $\mt$, for which the $\gamma_k$ are forced to lie in the fundamental alcove.

\noindent{\it Group actions}

The group $G^n$ acts on $\mg$ in the following way,
 defined  in (5.9) and 
(5.10) of [J1]:
An element $ \sigma$ in the first copy $G_1$ of $G$ acts on 
$\mg$ by sending
$$A_j   \mapsto \sigma A_j \sigma^{-1} ~{\rm for}~ j = 1, \dots, g $$
$$B_j   \mapsto \sigma B_j \sigma^{-1} ~{\rm for}~ j = 1, \dots, g $$
$$ D_j \mapsto \sigma D_j  \eqno (2.9)$$
$$ \gamma_1 \mapsto \sigma \gamma_1 \sigma^{-1} $$
$$ \gamma_j \mapsto \gamma_j ~{\rm for}~ j = 2, \dots, n  $$

while an element $\sigma$ in the $l$-th copy $G_l$ of $G$ 
(for $l = 2, \dots, n $) acts by 
$$ A_j \mapsto A_j, ~~~ B_j \mapsto B_j  $$
$$ D_j \mapsto D_j ~{\rm for}~ j \ne l $$
$$ D_l \mapsto D_l \sigma^{-1} \eqno(2.10) $$%{(*groupactmod)}
$$\gamma_j \mapsto \gamma_j ~{\rm for}~ j \ne l, $$
$$\gamma_l \mapsto \sigma \gamma_l \sigma^{-1} .$$
The action on $\mG$ is defined in the same way (replacing
the $\gamma_j$ by the $C_j$).
These actions  restrict to define  actions of $T$ on 
$\mt$ and $\mT$.

In terms of flat connections, the group actions on $\mg$ and $\mt$ result from an action of
the space of gauge transformations modulo those gauge transformations
which are trivial on the boundary: in other words the 
action results from changing the trivialization of the bundle
near the boundary components.

\noindent{\it Smoothness}

{\sc Proposition (2.11):} \tensl The spaces $\mG$ are smooth, and isomorphic
to $G^{2(g+n-1)}$.\tenrm
 
{\sc Proof:} The relation $F=1$ of (2.6) allows us to write $C_1$
as a function of the other variables $A_i, B_i, C_j,D_j$.

{\sc Proposition (2.12):} \tensl If 
$G = SU(2)$, the spaces $\mg$ are smooth near any point 
$(\{A_j, B_j\})$, $\{ \gamma_j, D_j\} )$ where 
either at least one of the $\gamma_j$ is such that  the exponential 
map is a local diffeomorphism in a neighbourhood of 
$\gamma_j$, (this is the case when $\gamma_j$ is regular, or $\gamma_j=0$),
  or the stabilizer
of the action of $G$  by conjugation at the point 
$(\{A_j, B_j\})$ is the centre $Z(G)$.\tenrm

{\sc Proof:} The proof of the 
previous proposition  shows that the differential of the
map $F \circ  E: G^{2g+n-1} \times \lieg^n \to G$  is surjective
if  the exponential map  is a local diffeomorphism in 
a neighbourhood of at least one of the $\gamma_j$.
Here, $F$ was defined by equation (2.3) and 
$$ E: (\{ A_j, B_j, \gamma_j, D_j \} ) 
\mapsto (\{ A_j, B_j, \exp \gamma_j, D_j \}). $$
If $\gamma_j$ is not regular then $\exp \gamma_j = \pm 1$; thus, 
if none of the $\gamma_j$ are regular, 
 all the $C_j$ are in the 
centre $Z(G)$ and so 
$C_1 \prod_{k=2}^n D_k C_k D_k^{-1} \in Z(G). $ Also, 
the image of the differential of the map 
$$R: G^{2g} \to G $$
given by 
$$ R(A_1, \dots, A_g, B_1, \dots, B_g) = 
\prod_{j = 1}^g A_j B_j A_j^{-1} B_j^{-1} $$ 
is the orthocomplement of the Lie algebra of the stabilizer of   
$(A_1, \dots, A_g, B_1, \dots, B_g) $ under the 
conjugation action (see [G2], Proposition 3.7).
If the stabilizer at the 
point $\{ A_j, B_j \} $ is 
$Z(G)$,  Goldman's result implies the differential 
of $R$ is surjective.  This completes the proof.

The above proof extends to give the following result.

{\sc Proposition (2.13):} \tensl If $G = SU(2)$ and 
any of the $C_j$ is regular for $j \ge 1$, 
or  the stabilizer of the conjugation action 
of $G$ at the point $\{ A_j, B_j \} $ $ \in G^{2g}$ is 
$Z(G)$, then $M^\liet$ is smooth at the point
$\{ A_j, B_j, \gamma_j, D_j \} $ and
$M^T$ is smooth at the point
$\{ A_j, B_j, C_j, D_j \}$. \tenrm

\noindent {\it Symplectic structures and imploded cross-sections}

In [J1] and [J2]  a closed 
2-form on $\mg$ was  constructed, for which the map 
$$\mu_k : (\{ A_j, B_j\} , \{\gamma_j,D_j\} ) \mapsto -\gamma_k $$
is the moment map for the action of the $k$-th copy of $G$. Similarly, 
for $\mt$, there is again a two-form for which $\gamma_k$ is the moment
map for the $T$-action. This form is invariant under $(W\ltimes\bbz^r)^n$,
and so descends to  $\mT$:

{\sc Proposition (2.14):}\tensl The space $\mg$ is equipped with a closed 
2-form, whose restriction to $\cap_{j = 1}^g \mu_j^{-1} (B)$ is
nondegenerate. Here, $B$ is the subset of $\lieg$ where the exponential
map is locally a diffeomorphism. 

The space $\mt$ is equipped with a a closed 
2-form, whose restriction to $\cap_{j = 1}^g \mu_j^{-1} (Reg)$ is
nondegenerate. Here $Reg$ is the subset of $\liet$ consisting of elements which exponentiate to 
regular elements of the torus; it consists of all translates under 
$(W\ltimes\bbz^r)^n$ of the interior of the fundamental alcove.

Quotienting, $\mT$ in turn is equipped with   a closed 
2-form, whose restriction to $\cap_{j = 1}^g \mu_j^{-1} (\Delta^0)$ is
nondegenerate. \tenrm

{\sc Proof:} See [J1], and Huebschmann [H] (Theorem 2.15). 
Huebschmann in fact
proves this result for a smooth  extended moduli space with a
Hamiltonian  action of 
$H \times G^n$, where $H \cong G$ and 
 the action of  $H$
is free (and hence $0$ is a regular value for its moment map): 
the appropriate subset of our
space $\mg$ is recovered by taking the symplectic quotient by the 
action of 
$H$.

The space $\mT$ is  the cross-section 
$\cap_{j = 1}^n (\mu_j)^{-1} (\Delta)$ of $\mg$.
 We shall need to construct auxiliary spaces
({\it imploded cross-sections}) where the degenerate directions of the
two-form have been collapsed. The definition of these spaces is
due to Sjamaar [Sj]. 

{\sc Definition (2.16):}\tensl  Let $M$ be a symplectic manifold equipped with the 
Hamiltonian action of a compact Lie group $G$ with maximal torus 
$T$ and 
moment map $\Phi: M \to \lieg$.  
Let $\sigma$ index the faces of the fundamental Weyl chamber.
Let $G_\sigma$ be the stabilizer of $\sigma$ under the adjoint action 
on $\lieg$. 
The imploded cross-section of $M$ is then 

$$ M_{\rm impl} = \coprod_\sigma {\Phi^{-1}(\sigma)\over 
[G_\sigma, G_\sigma]} $$
where $[G_\sigma, G_\sigma]$ is the commutator subgroup of 
$G_\sigma$. \tenrm

The imploded cross-section admits an action of $T$.

{\sc Theorem (2.17) }  [Sj] : \tensl The space $M_{\rm impl}$
inherits from $M$ the structure of a stratified symplectic space 
and a Hamiltonian $T$-action. The symplectic strata are 
the quotients $\Phi^{-1}(\sigma)/[G_\sigma, G_\sigma] $.\tenrm

A special case of the imploded cross-section is the imploded 
cross-section of the cotangent 
bundle $T^* G \cong G \times \lieg$. 
The cotangent bundle
admits an action of $G \times G$ as follows,
 corresponding to left and right multiplication.
If $(g_1, 1) $ denotes a point in the first copy of $G$ in 
$G \times G$, it sends a point 
$(h, X_{h})$ (where $X_{h}$ denotes the value at 
$h$ of the left invariant vector field generated by 
$X \in \lieg$) to 
$(g_1 h, X_{g_1 h})$. Likewise, the element 
$(1, g_2)$ of $G \times G$ sends the point
$(h, X_{h})$ to $(h g_2, ({\rm Ad} (g_2^{-1}) X)_{h g_2} ). $
The imploded cross-section of $T^* G$ (corresponding
to the action of the first copy of $G$)  is 
$$ \Epsilon (T^* G) = \coprod_\sigma {G\over [G_\sigma, G_\sigma]} \times \sigma, $$
where $\sigma$ ranges over the (interiors
of the) faces of the fundamental Weyl chamber, and $G_\sigma$ is the stabilizer
of $\sigma$ under the conjugation action. 
The imploded cross section of 
$T^*G$ inherits an action of $T \times G$.

When $G=SU(2)$, one obtains, for example, that $ \Epsilon (T^* G) = \bbc^2$,
with the standard two-form. In general, one has [Sj] that
the implosion $\Epsilon (T^* G)$ is only smooth over the locus where
$[G_\sigma,G_\sigma]$ is a product of $SU(2)$s.  

One finds [Sj] that the imploded cross section of any Hamiltonian
$G$-manifold $M$ is given by 
$$M_{\rm impl} = \left ( M \times \Epsilon (T^* G) \right) // G \eqno (2.18)$$
$$ = \{ (m,y) \in M \times \Epsilon (T^* G): \Phi(m) = \Phi_G(y)\}/G. $$
In this sense $\Epsilon (T^* G) $ may
be thought of as a universal imploded 
cross-section.
It gives the implosion of $M$ as 
the symplectic quotient of $M \times \Epsilon (T^* G)$ under the 
diagonal action of $G$ (where 
$G$ acts on $\Epsilon(T^*G)$ via  the action of $T \times G$, and 
we use minus the usual symplectic structure on $\Epsilon (T^* G)$): 
here we have written the moment map for the 
action of $G$ on $M $ as $\Phi$ and that for the action on 
$\Epsilon (T^*G)$ as $\Phi_G$.

Note that this implies that when $G=SU(2)$, the implosion of $M$ is
smooth over the locus where $G$ acts freely, or with a constant
stabiliser over an open set.

We may construct the imploded cross-section of 
$\mg$ with respect to the action of $SU(2)^n$ as 
$$ (\mg)_{\rm impl} =\coprod_{\sigma_1, \dots, \sigma_n}  
{ { \cap_{i = 1}^n (\mu_i)^{-1} (\sigma_i ) }
\over {[G_{\sigma_1}, G_{\sigma_1} ] 
\times \dots \times [G_{\sigma_n}, G_{\sigma_n} ]} }. \eqno (2.19)$$
This object may also be exhibited as the symplectic quotient of 
$\mg \times  \Epsilon (T^* G) \times \dots 
\times \Epsilon (T^* G)$  under the diagonal action of 
$G^n$.

\noindent{\it Quasi-Hamiltonian $G$-spaces}

We note that we have not discussed any symplectic structure for 
$\mG$; indeed, the appropriate  structure on $\mG$ is that of 
a {\it quasi-Hamiltonian $G$-space}. Following 
  Alekseev, Malkin and Meinrenken [AMM], we will give a definition
of such spaces, which are 
   equipped 
with group actions and group-valued moment maps.

{\sc Definition (2.20)} {[AMM]} \tensl A quasi-Hamiltonian $G$-space
is a manifold $M$ with a $G$-action together with an 
invariant 2-form $\omega$ and an equivariant map 
$\Phi: M \to G$ (the $G$-valued 
moment map) for which 

1. $$d \omega = - \Phi^* \chi$$
where $ \chi$ is the 3-form given in terms of the left-invariant
Maurer-Cartan form $\theta \in \Omega^1 (G) \otimes \lieg$ 
and the inner product $(\cdot, \cdot)$ on $\lieg$ 
by 
$$ \chi = (1/12) (\theta, [\theta, \theta]). $$ 

2. The map $\Phi$ satisfies
$$ \iota ({v_\xi}) \omega = (1/2) \Phi^*(\theta + \bar{\theta}, \xi) $$
where $\bar{\theta}$ is the right-invariant 
Maurer-Cartan form,  $( \cdot, \cdot)$ is the bi-invariant
inner product on $\lieg$ and $v_\xi$ is the vector field corresponding to $\xi$. 

3. For $x \in M$ the kernel of $\omega_x$ is given by 
$${\rm Ker} (\omega_x) = \{ v_\xi, \xi \in {\rm Ker} 
({\rm Ad} \Phi_x + 1) \} . $$\tenrm

{\sc Remarks}: (1) We note that if $G$ is Abelian, then $M$ is symplectic; in general, 
of course, it is not. 

(2) If a space is a quasi-Hamiltonian $G$-space, it is not necessarily
a quasi-Hamiltonian $H$-space for a subgroup $H$ of $G$.

An example of a quasi-Hamiltonian space is the double 
$D(G) = G \times G$, which is a quasi-Hamiltonian $G \times G$-space
with the action 
$$(g_1, g_2): (a,b) \mapsto (g_1 a g_2^{-1}, 
g_2 b g_1^{-1} ) \eqno(2.21)$$%{(*groupact)}
and the $G \times G$-valued moment map 
$\Phi: (a,b) \mapsto (ab, a^{-1} b^{-1}). $
(We refer to $G \times G$ as $G_1 \times G_2$, where $G_j = G$.)
The two-form $\omega$ is given by 
$$ \omega = (1/2) (\pi_1^* \theta, \pi_2^* \bar{\theta}) +
(1/2) (\pi_1^* \bar{\theta}, \pi_2^* \theta) . \eqno (2.22)$$
It is shown in [AMM] (Proposition 3.2) that 
$D(G)$ is a quasi-Hamiltonian $G\times G$-space.

It will be convenient to introduce an alternative
system of coordinates on $D(G)$: these are given by 
$(u,v)$ where $u = a$ and $v = ba$. 
In these coordinates the action of $(g_1, g_2) \in G \times G$ is 
$$ (g_1, g_2): (u,v) \mapsto (g_1 u g_2^{-1}, {\rm Ad}_{g_2} v) 
\eqno(2.23)%{(*groupactalt)}
$$
and the moment maps are 
$$ \Phi_1 (u,v) = {\rm Ad}_u v,  \Phi_2(u,v) = v^{-1}. \eqno(2.24)$$%{(*momentmap)}.$$

A second example of a quasi-Hamiltonian $G$-space
is the space $M^G$: 
there is a 
 quasi-Hamiltonian structure on the spaces
$\mG$ and $\mT$, for which the maps 
$$\Phi_k: (\{ A_j, B_j\}, \{ C_j \}, \{ D_j\} ) \mapsto C_k^{-1}\eqno (2.25)$$ 
are $G$-valued moment maps for the action of the $k$-th copy of $G$.
(See [AMM], equation (39), Section 9.2.) 

One may define reduction for 
quasi-Hamiltonian $G$-spaces: it is shown in [AMM] (Theorem 5.1) that 
if  $M$ is a q-Hamiltonian $G_1 \times G_2$-space
and $f \in G_1$ is a regular value of the moment map $\Phi_1: M \to G_1$
then the reduced space 
$$M_f = (\Phi_1)^{-1}(f)/Z_f\eqno (2.26)$$
is equipped with the structure of a quasi-Hamiltonian $G_2$-space.
(Here, $Z_f$ is the stabilizer of $f$ under the action of $G_1$
by conjugation.)

There is also a cross-section theorem for quasi-Hamiltonian
$G$-spaces. If $M$ is a quasi-Hamiltonian $G$-space, and 
$f \in G$, and $Z_f$ is the 
centralizer of $f$ in $G$, 
there is an open subset $U \subset Z_f$ for which the map 
$$G \times_{Z_f} U  \to G,\quad  [g, u] \mapsto gu\eqno (2.27)$$
is a diffeomorphism onto an open subset of $G$. It then follows 
([AMM], Proposition 7.1) that $\Phi^{-1} (U)$ is a smooth 
$Z_f$-invariant submanifold and is a quasi-Hamiltonian 
$Z_f$-space with the restriction of $\Phi$ as a moment map. Again,
if $Z_f$ is Abelian, the cross-section is symplectic.

\noindent{\it Imploded cross sections of quasi-Hamiltonian $G$-spaces}

Let $H,G$ be two Lie groups. We may formulate a  definition of 
imploded cross-sections of
quasi-Hamiltonian $H\times G$-spaces analogous to the definition for 
Hamiltonian $G$-spaces. Let the faces of $\exp \Delta \subset T \subset G$
(which are the same as the faces of $\Delta$)
be denoted by  $\sigma$,
 and let $G_\sigma$ be the stabilizer of the action of 
$G$ by conjugation at a point in the interior of $\sigma$. Then we define the 
imploded cross-section of a quasi-Hamiltonian $G$-space $M$ as 
$$M_{\rm impl} = \coprod_\sigma 
{\Phi^{-1}(\sigma)\over [G_\sigma, G_\sigma]}. \eqno (2.28)$$

We note that the generic locus of this implosion, that is,
 the stratum corresponding to the interior
of the alcove, is a quasi-Hamiltonian $H\times T$-space. 
In general, one would have to go through the analysis in [Sj] to show that 
one has a suitable   structure of a stratified $H\times T$-space on the full $M_{\rm impl}$. 
For 
 the case $G = SU(2)$, things can be somewhat simplified: 
the possible $\sigma$  are 
$\sigma^0 = \{ {\rm diag}  (e^{i \gamma}, e^{ - i \gamma}), 0 < \gamma < \pi \} $
(with $G_\sigma = T$, and $[G_\sigma, G_\sigma]$ trivial) and 
$\sigma^+ = 1 $ and $\sigma^- = - 1 $ (with 
$G_\sigma = [G_\sigma, G_\sigma] =G$); here 
$\pm 1$ refer to $\pm 1$ times the identity matrix in $SU(2)$.
The stratum of $M_{\rm impl} $ corresponding
to each $\sigma $ is a quasi-Hamiltonian $H\times T$-space,
 since according to the reduction theorem,
Theorem 5.1 of [AMM], the reduced spaces
$\Phi^{-1} (1)/G$ and 
$\Phi^{-1} (-1)/G$ inherit the structure of a
quasi-Hamiltonian $H\times T$-space (with $T$ acting trivially), and 
according to the cross-section theorem, Proposition 7.1 of [AMM], 
the cross-section $\Phi^{-1} (\sigma^0   ) $ also
inherits the structure of a
quasi-Hamiltonian $H\times T$-space. To see that the strata fit together correctly, 
we again, as above, construct a universal implosion, not from $T^*G$, but 
from the double $G\times G$, for $G= SU(2)$.

\noindent {\sc Proposition (2.29)}: \tensl For $G = SU(2), T=S^1$, the implosion 
$D(G)_{\rm impl}$ is smooth, isomorphic to $S^4$, and is a quasi-Hamiltonian 
$SU(2)\times S^1$-space. \tenrm

We have $D(G)$ as a quasi-Hamiltonian $G\times G$-space 
 and  
(under the action of $G = G_2$ in equation (2.23)) we obtain   
the  cross-section $\Phi^{-1}(\sigma^0)$ as 
$\{ (u,v) \in G \times G:  v^{-1} \in  \exp \Delta \} .$ 
The action of $g_2 \in G_\sigma  $ 
on $G \times \sigma$ is by right multiplication by $g_2^{-1}$
on $G$ and the adjoint action on $\sigma$ (the latter being
trivial by definition of $G_\sigma$). 
By the above discussion, the strata
$\Phi^{-1} (\sigma)/[G_\sigma, G_\sigma] $  are
quasi-Hamiltonian $G\times T$-spaces. 
The stratum $\sigma^0$ has stabilizer $T$, and $[T, T] $ is
trivial, so that $\phi^{-1} (\sigma^0)/[G_\sigma, G_\sigma] = 
G \times \sigma^0 = S^3\times \sigma^0. $
On the other hand the strata $\sigma^+$ and $\sigma^-$ have stabilizers 
$G_\sigma = [G_\sigma, G_\sigma] = G$, so for these strata 
we have $\Phi^{-1} (\sigma)/[G_\sigma, G_\sigma] = G/G $ which is 
a point (denoted $\pm 1$). 
The imploded cross-section is thus 
$$ D(G)_{\rm impl} = \coprod_\sigma  {G\over 
[G_\sigma, G_\sigma ] }\times \sigma $$
$$ = (G \times \sigma^0)\ \coprod  \{ 1 \} \coprod \{ -1 \}. $$
We may identify this with $S^4$ by identifying $\sigma^+$ with the 
north pole $N$ and $\sigma^-$ with the south pole $S$ , and 
projecting the cylinder $G \times (0, \pi)\simeq S^3\times (-1,1)$ in the 
standard fashion to the sphere  $S^4 - \{ N,S\}$.
By the cross-section theorem, $S^4 - \{ N,S\}$ inherits a $G\times T$ quasi-Hamiltonian
structure; one checks by explicit computation that it extends analytically
to the whole of $S^4$, and that it is non-degenerate at the poles.
 (Of course, this form is not closed.)
 The residual $G\times T$
action on $G\times \sigma^0$ is given by 
$$(g,t)\cdot (u,v) = (gut^{-1}, v),\eqno(2.30)$$
and corresponding moment maps:
$$\Phi_{G}: (u,v)\mapsto (uvu^{-1}),\ \Phi_{T}: (u,v)\mapsto  v^{-1}\eqno (2.31) $$
This action extends to $S^4$, fixing the poles.

We can use this universal   implosion to give  an alternate definition
of the implosion of any
quasi-Hamiltonian $G$-space, for $G= SU(2)$.

{\sc Proposition (2.32):} \tensl Let $M$ be a quasi-Hamiltonian $H\times G$-space (where
$G = SU(2)$) and
let $D(G)_{\rm impl} $ be the imploded cross-section of $D(G)$. Then 
$$M_{\rm impl} = M\times D(G)_{\rm impl}// G\eqno (2.33)$$
$$ = \{(m, \xi) \in M \times D(G)_{\rm impl}: \Phi(m) = \Phi_G(\xi)\}  /G. $$
The space $M_{\rm impl}$ is a quasi-Hamiltonian $H\times T$-space, and it is smooth over
the locus of points $(m,\xi)\in M\times D(G)_{\rm impl}$ where the stabiliser of the
 $G$-action is trivial.\tenrm

{\sc Proof:} We want to identify the stratum of $M_{\rm impl}$ corresponding to 
$\sigma$ with the set of points in the quotient $M\times D(G)_{\rm impl}// G$ 
with $\Phi(m)= \Phi_G(\xi) $ lying in the orbit of $\sigma$ under $G$. If 
$(m,\xi), \xi = (u,v)$ is such a point, we can normalise so that
 $\Phi(m) = \Phi_G(\xi) = uvu^{-1}\in \sigma$, in other words,  
$\Phi^{-1}(G\sigma)\times_{G\sigma} \Phi_G^{-1}(G\sigma)/G\simeq
\Phi^{-1}(\sigma)\times_\sigma \Phi_G^{-1}(\sigma)/G_\sigma$. 
One then has, for $(u,v)\in \Phi_G^{-1}(\sigma)$, that $uvu^{-1} = v$ (since $v$
lies in $\sigma$), and so
$u$ lies in $G_\sigma$, so that $\Phi_G^{-1}(\sigma)=  
(G_\sigma \times \sigma)/[G_\sigma,G_\sigma]$.
One wants to compare $\Phi^{-1}(\sigma)\times_\sigma \Phi_G^{-1}(\sigma)/G_\sigma\simeq
\Phi^{-1}(\sigma)\times_\sigma ((G_\sigma \times \sigma)/[G_\sigma,G_\sigma])/G_\sigma$
with $\Phi^{-1}(\sigma)/[G_\sigma,G_\sigma]$. This equivalence is straightforward.
On the open stratum, corresponding to the interior of the fundamental alcove,
we note that our comparison is between 
$\Phi^{-1}(\sigma)\times_\sigma \Phi_G^{-1}(\sigma)/T$
and $\Phi^{-1}(\sigma)$.

The smoothness of the quotient follows from the 
smoothness of $D(G)_{\rm impl} = S^4$. As the product $M\times D(G)_{\rm impl}$
 is a quasi-Hamiltonian $H\times G\times T$-space,
the reduction is a Hamiltonian $H\times T$-space, by the reduction theorem of [AMM]. 

\noindent{\it Imploded cross-section of $\mG$} 

We would now like to apply the implosion construction above to $\mG$, which,
as we saw, is a quasi-Hamiltonian $G^n$-space, where the moment map associated to the 
action of $G$ at the i-th puncture is just $C_i^{-1}$. Performing the implosion for 
each puncture, we will  denote
the resulting $\mG_{\rm impl}$ by P; it is a Hamiltonian $T^n$ space.

For $j = 2, \dots, n$, 
we find that we can first  identify the elements  
$\{ (D_j, C_j) \}$ of (2.6) with elements of 
$D(G)$: in the coordinate system  $(u,v)$
on $D(G)$,  $D_j$ corresponds to 
$u$ and $C_j$ to $v$. (See [AMM], section 9.2.)
It is clear that 
the imploded cross-section of $\mG$ (with respect to the action of the
$j$-th copy of $G$, where $j = 2, 
\dots, n$) is 
$M^0 \cup M^+ \cup M^-$, where $ M^0$ consists of the points in 
$\mG $ for which $(C_j)^{-1} \in \exp \Delta^0$, while 
$M^{\pm} $ are the quotients by 
$G$ of the sets of points for which $C_j = \pm 1$.  
If $C_j = \pm 1$, the action of $g \in G_j$ sends $D_j$ to 
$ D_j g^{-1}$ and preserves the values of $A_j, B_j, C_j$. 
Notice that for 
$j \ge 2$ the  $j$-th copy of $G $ acts on $\mG$ only through
its action on 
the variables $C_j$ and $D_j$.

It follows that the imploded cross-section of $\mG$ (after imploding 
in turn the actions of  the $n-1$ copies of $G$ corresponding to $j = 2, \dots, n$) is 
$$\eqalign{ \{(A_j,  B_j ,   C_1, W_2, \dots, W_n ):& A_k, 
B_k , C_1 \in G, ~~W_j \in (D( G))_{\rm impl},\cr  
&\prod_{j =1}^g [A_j, B_j] C_1  \prod_{j =2}^n (\Phi_G(W_j)) = 1 
  \} \cr}. $$
Here, $\Phi_G: D(G)_{\rm impl} \to G$ is the map 
$ (D, C) \mapsto  D C D^{-1}. $ 

Finally,
to complete the description of the imploded cross-section, we 
need to take the 
cross-section for the first copy  of 
$G$ (denoted $G_1$ in the paragraph preceding (2.9)) and collapse it appropriately.
The imploded cross-section is the union of strata given by 
$${ \{ \{ A_k, B_k \}, C_1, W_2, \dots, W_n: 
\prod_{j =1}^g [A_j, B_j] C_1 \prod_{j =2}^n (\Phi_G(W_j)) = 1, 
~C_1 \in \sigma
\} \over [G_\sigma, G_\sigma]} .$$
Setting $W_1 = (D_1, C_1)$, with an arbitrary $D_1$, 
the imploded  cross-section becomes 
$$ P= M^G_{\rm impl} = \{ \{ A_k, B_k \}, W_1, \dots, W_n : 
\prod_{j =1}^g [A_j, B_j] \Phi_G(W_1)    \prod_{j =2}^n \Phi_G(W_j)  = 1
\} /G,\eqno (2.34) $$
where the $W_j \in D(G)_{\rm impl}$ and 
the action of $g \in G$  is by 
$$(A_k, B_k, W_j= (u_j, v_j))\mapsto (gA_kg^{-1},g B_kg^{-1}, (gu_j, v_j))\eqno (2.35)$$

One may also understand this in terms of an extended moduli space
${\tilde M}^G$ $\cong G^{2g+2n}$ equipped with an action of 
$G^{n+1} = G_0 \times \dots \times G_n, G_i\simeq SU(2)$ (this corresponds 
to putting in an extra puncture): our space 
$M^G$ is the quasi-Hamiltonian quotient of ${\tilde M}^G$ with respect
to the action of $G_0$. The space ${\tilde M}^G$ is defined by 
$${\tilde M}^G = 
\{ A_1, \dots, A_g, B_1, \dots, B_g, C_1, \dots, C_n,D_1, \dots, D_n \} $$
and the action of $\sigma  \in G_0$ is given by
$$A_j   \mapsto \sigma A_j \sigma^{-1}, 
~~B_j   \mapsto \sigma B_j \sigma^{-1} ~{\rm for}~ j = 1, \dots, g $$
$$ D_j \mapsto \sigma D_j $$
$$ C_j  \mapsto C_j  $$
while the actions of the other $G_j$ (for 
$j = 1, \dots, n$) are given by the formulas 
(2.10). 
The moment map for the action of $G_0$ is then given by 
$$\Phi_0: (A_1, \dots, A_g, B_1, \dots, B_g, C_1, \dots, 
C_n, D_1, \dots, D_n) \mapsto \prod_{j = 1}^g A_j B_j A_j^{-1} B_j^{-1} 
\prod_{k = 1}^n D_k C_k D_k^{-1}. $$
The space $M^G_{\rm impl}$ is obtained by imploding 
the actions of $G_1 \dots, G_n$ and finally taking the quasi-Hamiltonian 
quotient by the action of $G_0$. 

{\it Weighted frames and the relation to $\mT$}

One obtains the imploded cross-section by taking the inverse image under the moment
map of the fundamental alcove, and then quotienting by a group which depends on which 
face of the alcove one is mapping to. In the case at hand, this inverse image is 
$\mT$, and we have a projection $\mT\rightarrow P$, which 
is an  isomorphism  over the interior $(\Delta^0)^n$ of the alcove.

{\sc Proposition (2.37)} \tensl Over $(\Delta^0)^n$, $\mT\rightarrow P$ is a symplectomorphism.
\tenrm

{\sc Proof:} For both $\mT$ and $\mG$, the symplectic, or quasi-Hamiltonian
 structures are best defined in terms of their description on spaces of flat connections.
For $\mT$, one can define the relevant space of flat connections in various ways, but 
the most convenient is that given in [J1], where one considers the space ${\cal A}^T$
of connections which are flat and  which are, in a neighbourhood of the punctures, 
constant, of the form $a\ d\theta$, with $a\in \liet$ and $\theta$ and angular variable
 around the puncture. If $\alpha, \beta$ are two $\lieg$-valued 1-forms representing
the variation of such connections, one defines the symplectic form
as $$ \sigma(\alpha, \beta) = \int_\Sigma \alpha\wedge\beta.\eqno (2.38)$$
For $\mG$, following [AMM], section 9, one uses a larger space of connections ${\cal A}^G$, which has well defined values on the boundary circle
surrounding the puncture. One then has restriction maps to the boundary circles
$$R_i: {\cal A}^G\rightarrow L\lieg,\eqno (2.39)$$
where $L\lieg$ is the loop algebra of $\lieg$. One then defines the two-form on 
${\cal A}^G$ by 
$$\hat \sigma(\alpha, \beta)  = \sigma(\alpha, \beta)  + 
\sum_iR_i^*\omega, \eqno(2.40)$$
where $\omega$ is a form on $ L\lieg$, given in [AMM], equation (37). The key point for us is that 
restricting to ${\cal A}^T$, we have that the map  $R_i$ takes its values in a one-dimensional
space of constant loops, and so $R_i^*\omega$ is zero on this space. This then means
that if $I: {\cal A}^T\rightarrow {\cal A}^G$ is the inclusion, inducing the map 
$\mT\rightarrow P$, one has that $I^*\hat \sigma= \sigma$, giving us our symplectomorphism.

An element of $\mT$ is a flat connection, with   trivialisations at the punctures such
that the holonomy lies in the fundamental alcove. One obtains $P$ from $\mT$ by 
collapsing some subvarieties, and one could ask what class of
geometric objects our space $P$ corresponds to. The answer can be seen as follows.
Over the interior of the alcove, one has $\mT\simeq P$, and so one again has 
a framing at each puncture, as before. As one moves to the boundary of the alcove,
i.e. to a holonomy of $\pm 1$, one finds in $T$ that any trivialisation gives 
a holonomy in the fundamental alcove. In $P$, all of these trivialisations are
identified; the trivialisation is ``collapsed''. One can think of it in these terms:
an element $(h, u)\in S^3\times [0,1]$ gives a framing ($h$), and a holonomy class
$u$ in the fundamental alcove; as one moves $u$ to the boundary (to 0 or 1),
 one collapses the framing (which can be thought of as a unit vector in $\bbc^2$), 
scaling it by a ``weight'' to zero, and obtains $S^4$; $S^4$ is a space of pairs
(``weighted'' frame, holonomy angle),  where the weight depends on the holonomy angle 
and goes to zero as the holonomy angle
moves to 0 or 1. The space $P$ then appears as a space of flat connections with these weighted
frames at the punctures.

{\bf 3. Glueing, and connections on the trinion.}

{\it i) Glueing.}

Now suppose that we have a  punctured surface $\Sigma_0$, possibly disconnected,
 with 
punctures $p_{+},p_{-}, p_1,...p_n$ and corresponding moduli spaces $\mT_0, P_0$.
 Choosing parametrised  boundary curves $c_+$ around $p_+$, $c_-$ around
$p_-$, we can build a surface $\Sigma$ by identifying 
these curves to one curve $c_\pm$; let $\mT, P$ be the moduli spaces
corresponding to $\Sigma$. 
The group $T=S^1$ embeds antidiagonally into $S^1\times S^1$; this group acts naturally
on the framings at $p_+, p_-$ of elements of $\mT_0$, with moment map 
the difference $t_+-t_-$ of the holonomy angles. We can take the symplectic quotient
$$\mT_{ 0,red} = \mT_0//T.\eqno (3.1)$$
(The symplectic form is of course degenerate when $t_+= \pm 1$.)
Geometrically, this amounts to asking that the difference of the holonomy angles $t_+-t_-$
vanish; one can then use the trivialisations to glue together the flat connections and 
obtain a flat connection on $\Sigma$. Along a $T$-orbit, the glueings give the same
 flat connection, so that there is a well-defined map
$$\mT_{0,red}\rightarrow\mT.\eqno (3.2)$$
The fiber of this map is a point if $t_+\ne 0, 1$. If $t_+=0$ or $1$, the fiber is  $S^3/S^1$
if $n\ge 1$; the reason is that the holonomy then lies in the maximal torus 
for any trivialisation, but one is only quotienting out a torus' worth of such choices.
If $t_+ = 0$ or $1$, and $n= 0$, the fiber is the quotient of $S^3$ by $S^1$ and the action of the group
of automorphisms of the bundle. 

The action of the subgroup $1\times S^1$ of $S^1\times S^1$ commutes with that of $T$, and so descends
 to the quotient; using the map (3.2), one has an action of $S^1$ on the moduli
 $\mT$
which is only defined when $t_+ \ne 0,1$. These are the flows of Goldman [G1].

In parallel to this, one can also take the symplectic quotient on the space $P^0$, obtaining 
$$P_{0,red} = P_0//S^1\eqno (3.3)$$
This space is symplectic. The action of the subgroup $1\times S^1$ descends to
$P_{0,red}$. Also, the map $\mT_{0}\rightarrow P_0$ descends to 
$$\mT_{0,red} \rightarrow P_{0,red}\eqno (3.4)$$
which commutes with the $1\times S^1$-action.

This construction can obviously be iterated; in particular, if $\Sigma$ is a  surface
with  moduli space $\mT$,  we can cut $\Sigma$  along $k$ circles 
into pieces $\Sigma_1, ...\Sigma_\ell$.
For any such circle, defining, say, the boundary betwen 
$\Sigma_i$ and $\Sigma_j$,  we have an action of $S^1\times S^1$ on the trivialisations of 
the bundles at the marked points of $\Sigma_i$ and $\Sigma_j$. If we take the symplectic quotients
by the  anti-diagonal $S^1$'s, we have  maps
$$ (\mT_1\times \mT_2\times ...\times \mT_n)//(S^1)^k \rightarrow 
\mT\eqno (3.5)$$
$$ (\mT_1\times \mT_2\times ...\times  \mT_n)//(S^1)^k \rightarrow 
(P_1\times P_2\times...\times P_n)//(S^1)^k\eqno (3.6)$$
The map (3.6) factors through $\mT$, giving a map 
$$\Psi: \mT\rightarrow (P_1\times P_2\times...\times P_n)//(S^1)^k\eqno (3.7)$$
Over the generic locus of points where all the holonomies at the punctures are not central, 
this map is an isomorphism; over a connection $E$ in $\mT$ where some, say $\ell$, of the holonomies are central, 
the fiber is the quotient of $(S^3)^\ell$ by the action of the automorphisms of $E$.

All of these maps commute with the $1\times S^1$-actions defined above;
 on ${\cal M}$, these gave 
the Goldman flows, which are only defined when the holonomy is non-central;
 on the right, we have globally defined  Hamiltonian $S^1$-flows. 

One case which is of prime interest is when we have decomposed the surface with no punctures
into a
sum of $2g-2$ trinions. In this case, as we saw, we almost had on $\mT$ the 
structure of a toric variety. We will see that the right hand side of (3.7) will then simply be the toric variety 
corresponding to $\mT$. We first  want to understand the structure of $P$ when 
$\Sigma$ is a trinion.

 \noindent{\it iii) Bundles on a trinion}

We now specialise to the case of $SU(2)$ bundles on the three-punctured sphere
or trinion. The first question is to analyse the space of representations
of the fundamental group. This group is the free group on two letters, but,
as in [JW], we want a description which is symmetric in the three punctures.
Let $k_1, k_2, k_3$ then be three generators, such  that $k_1 k_2 k_3=1$.
In terms of the preceding notation, we have
 $k_1 = c_1, k_2 = d_2c_2d_2^{-1} , k_3 = d_3c_3d_3^{-1}$.
We let $K_1, K_2, K_3$ denote the  images of $k_1, k_2, k_3$ under
a representation into 
$SU(2)$. If the $K_i$ are not equal to $\pm1$, they each determine in turn
a unique maximal torus, and indeed a unique translate in $\bbp_1=SU(2)/T$
of the fundamental chamber $\Delta= [0,\pi]$; these three translates are
 distinct for
an irreducible representation. More concretely, for 
each element, say $K_1$, there is a unique
point $\hat K_1$ in $S^2$ such that the action of conjugation by $K_1$ 
on $S^2 = SU(2)/S^1$
is given by an anti-clockwise
rotation $R_{2\gamma_1}^{\hat K_1}$ by an angle $2\gamma_1, \gamma_1\in (0,\pi)$
around $\hat K_1$. From the relation 
$$R_{2\gamma_1}^{\hat K_1}R_{2\gamma_2}^{\hat K_2}R_{2\gamma_3}^{\hat K_3 }=1,$$
one finds that the vertices $\hat K_1,\hat K_2, \hat K_3 $ form a spherical triangle with 
angles $\pi-\gamma_1, \pi-\gamma_2,\pi-\gamma_3$ respectively, and
 that $\hat K_3 $ can 
belong to only one  of the two hemispheres bounded
by the great circle determined by $\hat K_1$, $\hat K_2$ (there is an orientation).
 From these facts, or 
referring to [JW], one has the relations:
$$\eqalign{\gamma_1 +\gamma_2+\gamma_3\le& 2\pi\cr
\gamma_1 +\gamma_2-\gamma_3\ge& 0\cr
\gamma_1 -\gamma_2+\gamma_3\ge& 0\cr
-\gamma_1 +\gamma_2+\gamma_3\ge& 0\cr}\eqno(3.8)$$
This defines a tetrahedron $\Gamma$ in $\bbr^3$, and one has 

{\sc Lemma (3.9):} \tensl The equivalence classes of representations 
of the fundamental group of the trinion into $SU(2)$ are in bijective 
correspondence with elements of $\Gamma$.\tenrm 

The proof is given in [JW]. We note here that the action of conjugation on a 
representation is simply that of rotating the triangle in $S^2$. One can normalise
an irreducible representation so that $\hat K_1$ is the north pole,
(conjugating $K_1$ to $\Delta$), and then rotate around the north pole so that 
 $\hat K_2$ lies on a fixed meridian.

Several other remarks are in order.
\item {1.} The interior of $\Gamma$ corresponds to irreducible representations,
or again to non-degene\-ra\-te triangles in $S^2$.
\item {2.} The interior of the faces of the tetrahedron correspond
to reducible representations, with none of the $K_1,K_2,K_3 $ equal to $\pm1$.
In terms of triangles, they correspond to degenerate triangles, with 
either all three vertices at the same point, or two of the vertices at
one point and the other at the antipodal point.
\item {3.} The interiors of the edges correspond to representations at 
which one and only one of the $K_1,K_2,K_3 $ is central. There are six edges,
corresponding to $K_1=1, K_1=-1, K_2=1, K_2=-1, K_3 =1, K_3 =-1$ respectively.
\item{4.} The vertices correspond to representations for which 
all of the $K_1,K_2,K_3 $ are central. There are four vertices, corresponding 
to $(K_1,K_2,K_3 ) = (1,1,1), (1,-1,-1),$ $ (-1,1,-1),(-1,-1,1)$ respectively.

We now consider the spaces $\mT$, $P$ for the trinion, and group $SU(2)$.
We have 
$$ \mT=\{ (C_1, C_2, C_3, D_2, D_3)\in \Delta^3\times SU(2)^2| C_1D_2C_2D_2^{-1}
D_3C_3D_3^{-1}=1\}$$
In terms of representations, we have normalised $K_1$ to $C_1\in T$, and 
$K_2 = D_2C_2D_2^{-1}$, $K_3 =D_3C_3D_3^{-1}$. The fiber of this space over 
the representation space $R$, for an irreducible  element of $R$,
is given by 
$T^3/\pm (1,1,1)\simeq (S^1)^3$.
On the other hand, for the space $P$, one has fibers over 
the interior of $\Gamma$ of the form $S^1\times S^1\times S^1/\pm 1= (S^1)^3$. 
For representations in the interior of one of the faces (reducible 
representations with none of the $K_1,K_2,K_3 $ central), one has for 
each of $K_1,K_2,K_3 $ a frame in $S^1$, and one then quotients by the diagonal
$S^1$ to obtain as fiber an $S^1\times S^1$. When one goes to the edges,
one loses one of the framings as one of $K_1,K_2,K_3 $ is then central; the 
fiber is then a circle. Finally, at the vertices, there is no extra information
apart from $K_1,K_2,K_3 $ and so the fiber is a point.

In short, this is compatible with $P$ being the toric variety 
corresponding to $\Gamma$; this variety is $\bbp^3$. To see that this is the 
case, we first note 

{\sc Proposition (3.10)}: \tensl The map $\mu: P\rightarrow\Gamma$ which associates
to each representation the angles $\gamma_1, \gamma_2,\gamma_3$ is the moment map for
 the  action of $(S^1)^3$ on $P$.

The map  $\mu$ has a Lagrangian section, given as the 
fixed point set of an anti-symplectic involution. \tenrm

{\sc Proof:} The first statement is found in [J1], proposition 4.1. 
For the second part, it is easiest to give the involution on $\mT$ first.
To do this, we suppose that the three punctures on our Riemann sphere
all lie on the real equator $\bbr\bbp^1\subset \bbc\bbp^1$. The 
involution $z\mapsto\bar z$ induces a transformation on the 
space of representations, by pull-back. As the map reverses 
orientations, it is anti-symplectic.
 Choosing our base-points correctly,
we can give the action on $\mT$, up to conjugation, by 
$$ (C_1,C_2,C_3, D_2, D_3)\mapsto 
(C_1^{-1}, C_2^{-1}, C_3^{-1}, C_1D_2, D_3)$$
This does not leave the $C$'s in $\Delta$, so we compose with 
conjugation by an element $J$ representing the non-zero Weyl group 
element, which gives $JC_i^{-1}J^{-1} = C_i, i= 1,2,3$. 
The result of the composition is then 
$$ (C_1,C_2,C_3, D_2, D_3)\mapsto (C_1,C_2,C_3,JC_1 D_2J^{-1}, 
JD_3J^{-1})$$
The induced action on $K_1,K_2,K_3 $ is $(K_1,K_2,K_3 )\mapsto (K_1, K_1^{-1}JK_2^{-1}J^{-1}K_1, 
JK_3 ^{-1}J^{-1})$. Let us consider this action for the irreducible 
representations in the interior of $\Gamma$.
Choosing polar coordinates $(\theta, \phi)\in [0,\pi]\times [0, 2\pi]$
 on the sphere 
so that $J$ corresponds to the point $\hat J = (\pi/2, 0)$,
the elements $K_1,K_2,K_3 $ are fixed if $\hat K_3$ lies on the great circle $\phi =\pm \pi/2$,
and $\hat K_2 $ lies on the great circle  $\phi = \pm(\gamma_1+\pi/2)$.
Taking orientations into account, this gives two possibilities for 
$(K_1,K_2,K_3 )$ once the representation has been fixed. Once $ K_1,  K_2,
 K_3 $ are fixed, 
one finds that there are two possibilities each for $D_2, D_3$, giving 
eight possibilities in all in the fiber $(S^1)^3$ of $P$
over $\Gamma$. The real points in the  fiber
are an orbit  of $\bbz/2\times \bbz/2\times \bbz/2$
in $S^1\times S^1\times S^1$.  Over the interior of $\Gamma$, we can normalise
our choices so as to get a Lagrangian section. This section extends in ${\cal M}_\Delta$ to the 
interior of the faces of $\Gamma$. As  the involution preserves conjugacy classes,
one checks that the involution 
descends to $P$, and here, 
it  extends to all of $\Gamma$. 
 From all of this, we obtain:

{\sc Proposition (3.11)}:\ \tensl The variety $P$ is isomorphic to $\bbp^3(\bbc)$.
\tenrm

{\sc Proof}: This is in essence a special case of Delzant's theorem [D].
Indeed, $\Gamma$ is the moment polytope for $\bbp^3(\bbc)$ with the 
action of $(a,b,c)\in (\bbc^*)^3$ given in inhomogeneous 
coordinates by 
$(a,b, c)\cdot (x, y, z) = (abx, bcy, acz)$.
 We first note that the variety $P$ is smooth, 
except possibly at four points corresponding to the vertices.
This follows from the fact, referring to (2.34), that $P$ is defined as the quotient of 
a subvariety of $(D(G)_{\rm impl})^3= (S^4)^3$, by the action of $SU(2)$, and 
the action of the group is free away from these points. Leaving these four points aside
for the moment, both $P$ and $\bbp^3(\bbc)$
have Lagrangian sections for their respective moment maps to $\Gamma$, and we can identify 
these sections. This identifies our  section of $P$ above with the positive octant
in $\bbp^3(\bbr)\subset \bbp^3(\bbc)$. 
The group actions are however compatible, and so we can use the action 
to identify $P$ and $\bbp^3(\bbc)$ globally.

 To see that the quotient space $P$ is smooth
even at the four vertex points, we first 
map a neighbourhood of the origin in $\bbc^3$ to $(S^4)^3$, in a way that
is generically bijective on the quotient. Let us do this, for the
point corresponding to the trivial flat connection.
We first extend the real involution to $(S^4)^3$. For this, 
we write $S^4$ as the result of collapsing the ends of $\Delta \times SU(2)$.
The real involution on $(S^4)^3$ is then obtained from the involution on 
$(\Delta \times SU(2))^3$: 
$$I:  (C_1,C_2,C_3, D_1,  D_2, D_3)\mapsto (C_1,C_2,C_3,JD_1J^{-1}, JD_1C_1 D_1^{-1}D_2J^{-1}, 
JD_3J^{-1})$$
The action of $g\in SU(2)$  maps $D_i$ to $gD_i$, and fixes the $C_i$; the action of $h\in S^1$, 
on the $j$th factor,
is given by $(C_j, D_j)\mapsto (C_j, D_jh)$. Note that $I\circ g = JgJ^{-1}\circ I$. Recall also that 
in $S^4$, the $D_i$'s get collapsed if $C_i = \pm 1$. Now lift the real octant 
in $P$ to $(S^4)^3$ by setting  $D_1 =1$, and then extend this via the $(S^1)^3$ actions to a
map of $\bbc^3$. This gives the slice, and so expresses the $SU(2)$-invariant
functions on $(S^4)^3$ as a subring of the functions  on $\bbc^3$.

One can then build invariant coordinates on $P$: set
$C_j = exp(diag (i\theta_j, -i\theta_j))$, and let
$r_1= \theta_2+\theta_3- \theta_1,  r_2= \theta_1+\theta_3- \theta_2,
r_3= \theta_1+\theta_2- \theta_3$. One considers the functions
$$\eqalign{ V_1 &= r_1 C_1^{1/2}D_2\pmatrix{ 1\cr 0}
\wedge D_3\pmatrix{ 1\cr 0},\cr
V_2 &= -r_2 D_3\pmatrix{ 1\cr 0}
\wedge D_1\pmatrix{ 1\cr 0},\cr
V_3 &= r_1 C_1^{1/2}D_2\pmatrix{ 1\cr 0}
\wedge D_1\pmatrix{ 1\cr 0},}$$
with values in $\Lambda^2(\bbc^2) = \bbc$.
These functions are invariant, and provide a 
$(S^1)^3$-equivariant identification 
of $P$ with an open set of $\bbc^3$, showing that the 
subring is in fact the whole ring of functions on $\bbc^3$.

Let us summarise the situation of the Goldman flows for a Riemann surface $\Sigma$ of genus $g$,
with $n$ punctures $p_1,\dots, p_n$. The moduli space $\mT$ is then of dimension $6g+4n-6$.
 We can cut the Riemann surface along $3g-3+n$ disjoint circles, decomposing the
surface into trinions, with at most one of the original punctures in each 
trinion;
the ends are either the cut circles or  circles 
around the punctures, that is $3g-3+2n$ circles in all; let $D$ denote this decomposition. Quotienting by the antidiagonal
$S^1$-actions at each of the cut circles, there is then a map, as in (3.7) 
$$\Psi: \mT\rightarrow  P^D = (\bbp^3(\bbc)\times \bbp^3(\bbc)\times....\times\bbp^3(\bbc))//(S^1)^{3g-3+n} 
\eqno (3.12)$$
which is generically an isomorphism. Here there are $2g-2+n$ copies of $\bbp^3(\bbc)$ in 
the product, corresponding to the same number of trinions in the decomposition.
On the right hand side, one has an effective action 
of $(S^1)^{3g-3+2n}$, and so the variety is toric. Let $\Gamma\subset P^D$ denote the set of elements
in $P^D$ corresponding to  flat  connections whose holonomy is
 central along one of the $3g-3+n$ cut circles,
and $\Gamma' $ the locus corresponding to connections which are central along all of the
$n$ puncture circles. We can choose the decomposition $D$ in such a way that the 
flat connections which are reducible all lie in $\Psi^{-1}(\Gamma)$, for example
letting one of the cut circles be null-homologous; 
let us call these good decompositions.
One can reduce both $\mT$ and $P^D$ under the action of the $S^1$'s
associated to the punctures. On $\mT$, 
 the resulting $\mT_{\rm red}$,  
for values of the moment map (holonomy) $\mt\rightarrow exp(\Delta)$
 different from $\pm 1$ at each puncture, is a parabolic moduli space. From 
$P^D$, one obtains a variety $P^D_{\rm red}$ which is again toric,
and one has a map 
$$\Psi_{\rm red}: \mT_{\rm red}\rightarrow  P^D_{\rm red}.\eqno (3.13)$$

Finally, we note that we can restrict the holonomy at one of the punctures, say $p_1$, 
to be $-1$. This corresponds to reducing at the singular value $-1$ 
of the moment map; doing this for $P^D$, one obtains a space $P_{-}^{D}$.
If we then reduce at regular values of the moment map at the other punctures, we have a
space $P^{D}_{\rm -, red}$. One could do the same for $\mT$; it is preferable
however to use the space $\mT_{-} = \mT_{-1}$ of (2.8)  associated to the surface 
punctured at $p_2,...,p_n$. There is a map 
$$\Psi_{-}:\mT_{-} \rightarrow P_{-}^{D},\eqno (3.14)$$
which is an isomorphism away from the varieties $\Gamma_{-}, \Gamma_{-}'$ of elements 
whose holonomy is central at one of the cut circles, or all of puncture circles, respectively.
Reducing at the other punctures, at regular elements, one has spaces and
a map: 
$$\Psi_{\rm -, red}:\mT_{\rm -, red} \rightarrow P^{D}_{\rm -, red}.
\eqno (3.15)$$
Summarising, one has:

  {\sc Theorem}: (3.13) \tensl a) 
For $n>0$, the variety $\mT$ is smooth away from 
$\Psi^{-1}(\Gamma')$. For any $n$ and a good decomposition
 it is smooth away from $\Psi^{-1}(\Gamma\cap \Gamma')$.
The variety $P^D$ is smooth away from $\Gamma$. 

The map $\Psi$ of  (3.12),
 away from $\Psi^{-1}(\Gamma\cup \Gamma')$,
is a symplectomorphism  between $\mT$ and the toric variety $P^D$. 
Over an element  $E$ in $\Gamma\cup \Gamma'$ where some, say $j$, of the holonomies are central, 
the fiber is the quotient of $(S^3)^j$ by the action of the automorphisms of $E$.

The real codimension of $\Gamma$ in $P^D$ is six; in $\mT$,
 the codimension 
of $\Psi^{-1}(\Gamma)$ or of $\Psi^{-1}(\Gamma')$ is three. 

b) For $n>0$, reducing by the $(S^1)^n$'s acting on the framings at the punctures, 
at holonomies different from $\pm 1$, the resulting $\mT_{\rm red}$ is smooth,
and is symplectomorphic to $P^D$ away from the reduction of
 $\Psi^{-1}(\Gamma) $. $P^D_{\rm red}$ is smooth away from $\Gamma$.

The map $\Psi_{\rm red}$ is a symplectomorphism, away from the reduction
of $\Psi^{-1}(\Gamma)$.
 The real codimension of the reduction of $\Gamma$ in $P^D_{\rm red}$ is six, and that
of the reduction of $\Psi^{-1}(\Gamma)$ in $\mT_{\rm red}$ is three.

c) The spaces $\mT_{-}$ are smooth away from $\Psi_{-}^{-1}(\Gamma_{-}')$,
and the $ \mT_{\rm -,red}$ are smooth, again reducing at holonomies different from
$\pm 1$ at the punctures $p_2,...,p_n$. 
The map $\Psi_{-}$ of  (3.14),
 away from $(\Psi_{-})^{-1}(\Gamma_{-}\cup \Gamma_{-}')$,
is a symplectomorphism  between $\mT_{-}$ and the toric variety $P_{-}^D$.  
The map $\Psi_{\rm -,red}$ of (3.15) is a symplectomorphism, away from the
 reduction of $(\Psi_{-})^{-1}(\Gamma_{-})$.

The real codimension of $\Gamma_{-}$ in $P_{-}^D$ is six; in $\mT$,
 the codimension 
of $\Psi_{-}^{-1}(\Gamma_{-})$ or of $\Psi_{-}^{-1}(\Gamma_{-}')$ is three.

 The real codimension of the reduction of $\Gamma_{-}$ in $P^D_{\rm -, red}$ is six, and that
of the reduction of $\Psi_{-}^{-1}(\Gamma_{-})$ in $\mT_{\rm -, red}$ is three.
\bk 

\tenrm

{\bf 4. Moduli of framed parabolic sheaves}

If one fixes the conjugacy class of the 
holonomy at the punctures, there is a well-established correspondence,
due to Mehta and Seshadri, between the moduli spaces of representations of $\pi_1(\Sigma)$
into $SU(2)$ on one hand, and, on the other,
 rank-2 holomorphic bundles with trivial determinant and with a parabolic structure at the 
punctures (a choice of a line $\ell_p$ in the fiber above the puncture $p$, along with 
some weights). 
The eigenvalues $exp (\pi i\delta), exp(\pi i(1-\delta)),
\delta\in [0, 1]$ of the holonomy get translated into weights  of the 
parabolic structure, and so to different choices of polarization. In our case, we are dealing with a space in which all of the
possible holonomies appear: we will thus want a space which contains all of these 
parabolic moduli spaces. We add in the extra "conjugate" variables of (partial) 
framings at  
the punctures, and so 
consider a space of framed parabolic bundles, that is bundles equipped with 
a trivialisation of the lines $\ell_p$. As our bundles $E$ have an $SL(2,\bbc)$ structure, 
this is the same as a trivialisation of  $E_p/\ell_p$, and so can be thought of as a map
$E\rightarrow \bbc_p$, whose kernel at $p$ is the line $\ell_p$. Moduli of such 
pairs  (bundles, maps to a fixed sheaf) have been studied by Huybrechts and Lehn
[HL] and we want to adapt their work. Our moduli spaces will have an 
extra parameter, corresponding to an $SL(2,\bbc)$ structure on the bundle.

We first recapitulate and make some comments on the moduli with fixed holonomy.
It will turn out that for $\theta= 1$, our moduli space, if it is to be compact,
 will contain  sheaves with torsion;
we will want to adapt a notion of stability  of [HL] for pairs consisting of 
a coherent sheaf and a map of the sheaf  into a fixed sheaf,
 which extends in our case  
the stability used for parabolic bundles.
Let $p_1,..,p_n$ denote the punctures and let $S$ be the sheaf $\oplus_i\bbc_{p_i}$. 
A (parabolic) framing of a sheaf $E$ is then a map of sheaves $\alpha: E\rightarrow S$, 
which we can write as $\alpha_i:E\rightarrow\bbc_{p_i}$. Let $\gamma_1,..,\gamma_n
\in [0,1]$ be a set of weights. For a subsheaf $F$ of $E$, we set 
$\mu_i(F) = 1$ if $F$ lies  
in the kernel of $\alpha_i$ at $p_i$, and $\mu_i = 0$ otherwise. The cases considered
in [HL] had $\alpha_i\neq 0$; we extend here to $\alpha_i = 0$. Let $\sigma_i(F) = 
1/rk(E)$ if $F=E$, and $0$ otherwise.
 We will say that
a pair $(E,\alpha)$ is (semi) stable if  for all non-trivial coherent subsheaves
 $F$ (including $E$),
$$rk(E)deg(F)<(\le) rk (F)(deg(E)-\sum_i\gamma_i) +
 rk (E)\sum_i(1-\mu_i(F) + \sigma_i(F))\gamma_i\eqno (4.1)$$
In our case, $rk(E)=2, deg(E) =0$, and so this becomes
$$2deg(F)<(\le) \sum_i\gamma_i(2-rk(F)-2\mu_i(F) +\sigma_i(F))\eqno (4.2)$$

\ni {\sc Lemma (4.3)} \tensl If $(E,\alpha)$ is semi-stable, then 

\item{i)} The kernel of $\alpha$ is torsion free, and the torsion subsheaf of $E$ is 
concentrated over the $p_i$, and is either $0$ or $\bbc$ at each $p_i$,
\item {ii)} If $\gamma_i> 0$, then the map $\alpha$ is non-zero at $p_i$,
\item{iii)} If $\gamma_i<1$, then $E$ is torsion free at $p_i$,
\item{iv)}  For $\gamma_i\in (0,1)$, one has a parabolic structure at $p_i$, 
and if all the  weights $\gamma_i$ lie in $(0,1)$, the stability condition (4.2) is identical to that
for parabolic bundles with weights $(1-\gamma_i)/2, (1+\gamma_i)/2$,
\item{v)} Let $(E,\alpha)$ be a semistable pair, which is locally free at $p_i$,
with a non-zero $\alpha$. For $\gamma_i = 0$, there is a family $(E_t, \alpha_t),t\in \bbc$ 
of semistable pairs such that $(E_t, \alpha_t)\simeq (E,\alpha), t\ne 0$, 
and $\alpha_0=0$.  For $\gamma_i = 1$,
there is a family $(E_t, \alpha_t),t\in \bbc$ 
of semistable pairs such that $(E_t, \alpha_t)\simeq (E,\alpha), t\ne 0$,
and $E_0$ has torsion at $p_i$.\tenrm

{\sc Proof}: For i), the torsion subsheaf of $ker(\alpha)$ has rank 0,   so 
its degree must be zero, and so it must be zero [HL]. For ii), one simply considers
the stability condition for  $E$ as a 
subsheaf of $E$. For  iii), one considers the torsion subsheaf of $E$ at $p_i$.
The existence of a parabolic structure in iv) follows from ii) and iii); 
the stability condition for parabolic bundles is 
$deg(G) + \sum_i(\mu_i(1+\gamma_i)/2 + (1-\mu_i) (1-\gamma_i)/2) <(\le) n/2$,
which is equivalent to (4.2) for line bundles, and the stability conditions for other 
subsheaves is automatically satisfied. 
For v), when $\gamma_i=0$, one simply sets $(E_t,\alpha_t)= (E,t\alpha)$.
For $\gamma_i = 1$, we choose a local coordinate $z$ near $p_i$, with $p_i$ 
corresponding to $z=0$. Over, say, $|z|<1$, we then choose a local isomorphism $E\simeq {\cal O}\oplus{\cal O}$,
such that $\alpha$ is represented by the vector $(1,0)$. We can write 
${\cal O}\oplus{\cal O}$ as a quotient, for $t\ne 0$:
$$\matrix{ &\pmatrix{-z&t&0}&&A=\pmatrix{t&z&0\cr0&0&1}& \cr {\cal O}&
\longrightarrow&{\cal O}\oplus{\cal O}
\oplus{\cal O}&\longrightarrow&{\cal O}\oplus{\cal O}}\eqno (4.4)$$
The map $t^{-1}\alpha = (t^{-1},0)$ can be represented by the map $(1,0,0)$ on 
${\cal O}\oplus{\cal O}\oplus{\cal O}$ at $z=0$. For $t=0$, the quotient in 
(4.4) is not ${\cal O}\oplus{\cal O}$, but $\bbc_{p_i}\oplus{\cal O}\oplus{\cal O}$.
One builds the family $E_t$ by using the matrix $A$
in (4.4) to map ${\cal O}\oplus{\cal O}\oplus{\cal O}$ to $E$ over $z\ne 0$; as this annihilates 
the kernel ${\cal O}$,   it  glues the quotient sheaf. The map $\alpha_t$ is then represented by 
$(1,0,0)$. One checks that $E_0$ is still semistable, using  a natural map $E_0/{\rm torsion}\rightarrow
E_t, t\ne 0$; this maps a subsheaf $L$ of $E_0$ to a subsheaf $\tilde L$ of $E_t$, killing the torsion,
but leaving the difference of the degree  and $\mu_i$ unchanged.

{\sc Remark}: As a consequence of v), any point in the moduli space will have a representative 
with torsion at $p_i$ for $\gamma_i=1$. At such points one can write the 
sheaf as a sum of its torsion piece and the kernel of $\alpha_i$. If one quotients out the torsion at the $k$
points with $\gamma_i$ = 1, one obtains a bundle with $c_1(E) = - k$, and the 
parabolic structure has disappeared.  
 Similarily, for $\gamma_i=0$, we can choose a representative
with $\alpha_i=0$, so again the parabolic structure has disappeared. If one does
this for both the points for which $\gamma_i=0$, and for which $\gamma_i=1$, one obtains
a bundle of degree $-k$ with parabolic structure at the $k_i$ remaining points,
and one checks that our stability condition becomes the stability condition
for such bundles.

We now recall some of  the construction of Huybrechts and Lehn, modified
 to correspond to our case, and in particular taking into account 
the existence of several punctures; the weights, initially, are in $(0,1)$,
so that one is dealing with vector bundles, and the $\alpha_i$ are all non-zero.
  The construction also gives the moduli of 
parabolic bundles; see also [Th], [Bh], and [Gi].
 As with any moduli space involving vector bundles, the starting point is the 
Grothendieck $Quot$-scheme parametrising quotients of the trivial bundle
$${\cal O}^{\oplus N}\rightarrow \tilde E\eqno (4.5)$$
(We twist $E$ by a fixed line bundle $L$  so that the resulting $\tilde E$
 is generated by global sections and
so that  the maps (4.5) induce isomorphisms on global sections; as we were dealing 
originally with $SL(2,\bbc)$-bundles, our bundles now have an induced isomorphism
 $\gamma: \Lambda^2(\tilde E)\simeq L^2$.)

In other words, we are not only considering $(E,\alpha_i)$, but have added in
a basis of sections of $E$; one has a good scheme parametrising such objects,
 but we now must quotient by the group $Gl(N,\bbc)$.

It is   convenient, again following [HL, Th, Bh, Gi], to transform the datum encoding $E$
somewhat. One takes the induced map on  second exterior powers, and then the map that this 
induces on sections; using the map $\gamma:\Lambda^2(E)\rightarrow L^2$
 given by the $Sl(2,\bbc)$ structure,  one gets an element $\hat \beta$ of 
$$V_1 = Hom(H^0(\Lambda^2{\cal O}^{\oplus N}), H^0(L^2)).
\eqno(4.6)$$
This map is  a finite injective morphism ([Th], section 7).

Next, one  adds in the parabolic data: each non-zero framing, acting on the global sections,
 yields an  element $\hat \alpha_i$ in copy $V_{2,i}$ of  
$$ V_2  = H^0({\cal O}^{\oplus N})^* ,\eqno (4.7)$$
for each puncture $p_i$.
 
If, as in the case of parabolic bundles, one is only interested in the 
$\alpha_i$ up to (independent) scale, one then obtains a point 
representing the equivalence class of ($E$ with its parabolic structure)
in a closed subvariety $X$ in  
$$Z_1\times Z_{2,1}\times Z_{2,2}\times\cdots\times Z_{2,n},\eqno (4.8)$$
where $Z_1=\bbp(V_1), Z_{2,i} =\bbp(V_{2,i})$, 
under the natural action of $Sl(N,\bbc)$ on these spaces. The space $X$ is
independent of the choice of weights; what varies is the polarisation, that is 
the choice of line bundles on which the action linearises.
Let 
$\pi_1,\pi_{2,i}$ denote the projections onto $Z_1, Z_{2,i}$ respectively. 
Set $L_{0} = \pi_1^*({\cal O}(N)), L_{1,i} = 
\pi_1^*({\cal O}(N-1))\otimes \pi_{2,i}^*({\cal O}(2 ))$; then the linearisation
corresponding to the (rational) weights $\gamma_i$ is $(L_{0})^{s_{0}}\otimes (
\otimes_i ( L_{1,i})^{s_{1,i}})$,
where  $s_{0}(\gamma_i)= s_{1,i}(1-\gamma_i)$. 

We now pass to the problem of fitting all of these moduli together, and adding in the 
cases $\gamma_i = 0,1$. 
 Fortunately, 
there is a technique due to Thaddeus [Th] which allows us to consider 
all of these quotients at once.
One builds over the space $X$ the $(\bbp^1)^n$-bundle
$$Y= \bbp(L_{0 }\oplus L_{1,1})\oplus \bbp(L_{0 }\oplus L_{1,2})
\oplus\cdots
\bbp(L_{0 }\oplus L_{1,n}),\eqno (4.9)$$ endowing it with its natural polarisation
${\cal O}(1,1,...,1)$. This space has embedded in it all the stable points for the 
various $\gamma_i$, and one can obtain the moduli spaces for each $\gamma_i$
by taking an appropriate $Sl(N,\bbc)\times (\bbc^*)^n$-quotient; the $\gamma_i$
now appear as the parameters for the appropriate $(\bbc^*)^n$ polarization [Th].
We will be interested, however, in the $Sl(N,\bbc)$ quotient of the 
whole ``master'' space $Y$; it will turn out that this is not quite the moduli
space which we want, but it is almost the correct one. One must modify the definition 
when $\gamma_i = 0$.

Indeed, we consider the space of quadruples
 $(E,\alpha_i,< \alpha_i>, \gamma)$, where $E$ is a rank 2 sheaf, 
$\alpha_i$ is the map $E\rightarrow \oplus_i\bbc_{p_i}$, 
$<\alpha_i>$ is a subspace
of $E|_{p_i}$ which is the kernel of $\alpha$ when $\alpha$ is non-zero, 
one-dimensional 
when $E$ has no torsion at $p_i$. This adds a parabolic structure when $\gamma_i=0$;
one has a projective class $<\hat \alpha_i>$ 
 of the corresponding map $\hat \alpha_i$ even when $\alpha_i$ vanishes. Our quadruples
then give elements $<\hat \beta>,<\hat \alpha_1>,...,<\alpha_n>$
in $X\subset Z_1\times Z_{2,1}\times Z_{2,2}\times\cdots\times Z_{2,n}$.

There are additional questions to be dealt with when $\gamma_i$ goes to 1; 
the bundles are then not necessarily locally free, and so can acquire torsion.
 We examine  the $Sl(2,\bbc)$-structure of such sheaves $E$.

{\sc Lemma (4.10)}: \tensl Let $E_t, t\in \bbc$ be a coherent family of rank 2  sheaves over the curve $\Sigma$, 
with $E_t$ locally free at $p$ for $t\neq 0$, and $E_0$ with $\bbc_p$ as  torsion subsheaf
near $p$. Let $\phi_t\in H^0(\Sigma, \Lambda^2(E)^*)$ be a family of $SL(2,\bbc)$-structures
on $E_t$. Then $\phi_0$ vanishes at $p$.\tenrm

{\sc Proof}: The question is local. If $z$ is a coordinate on $\Sigma$ with $z=0$ corresponding 
to $p$, one can obtain $E_t$ locally (changing the parameter $t$ if necessary),
as the cokernel of  
$${\cal O}{\buildrel {(0,t^k,z)}\over\longrightarrow}{\cal O}\oplus
{\cal O}\oplus{\cal O},$$
for some integer $k$. The forms $\phi_t$ are then given as multiples of  $e_1^*\wedge (-ze_2^*+t^ke_3^*)$, which vanishes
at $z=t=0$.

What this tells us is that there is up to scale a fiducial two-form $\phi$, which is non-zero 
away from the torsion points and vanishes at the torsion points. The choice of an $Sl(2,\bbc)$
structure is the choice of the (non-zero) scale for such a form.

One must see how our quadruples correspond to orbits  in $Y$, 
living over the orbits 
of $\hat E_X= (<\hat \beta>,<\hat \alpha_1>,...,<\hat \alpha_n>)$ in $X$. We note that as we
are dealing with projective bundles, we can tensor with line bundles without 
changing the structure, so that $Y$ can be written as the $(\bbp^1)^n$-bundle
$$\eqalign{ Y= \bbp(\pi_1^*({\cal O}(-1))\oplus \pi_{2,1}^*({\cal O}(-2)))\oplus &
\bbp(\pi_1^*({\cal O}(-1))\oplus \pi_{2,2}^*({\cal O}(-2)))
\oplus\cdots\cr
&\bbp(\pi_1^*({\cal O}(-1))\oplus \pi_{2,n}^*({\cal O}(-2))),}\eqno (4.11 )$$
 so that one is considering  projectivisations of the 
sum of a tautological bundle and a second power of a tautological bundle.
We then have a natural lift of $\hat E_X$ to $Y$, given by 
$$\hat E_X' =( (\hat \beta ,\hat \alpha_1^{2}),...,(\hat \beta,\hat \alpha_n^{2}))$$
This must be extended to the cases when the sheaf ceases to be locally free
and  acquires torsion at the $p_i$.
Indeed, when the bundle acquires torsion, at $p_i$, say, one finds that by rescaling 
the torsion piece of the bundle we can modify the scale of $\hat \alpha_i^2$ to 
$c\hat \alpha_i^2$, say. These should be represented in the same orbit, and 
if the lift
$\hat \beta$ were non zero, this would not be  the case. One would want the $i-th$ pair 
in $\hat E_X'$ to be 
of the form $(0, \hat \alpha_i^2)$ which would then absorb the scale. This can be achieved
as follows. Recall that $\hat \beta $ is defined as the composition.
$$\Lambda^2(\bbc^N){\buildrel \sigma\over 
\rightarrow } H^0(\Sigma,\Lambda^2(E)){\buildrel \xi\over 
\rightarrow }H^0(\Sigma,L^2).$$
The second map $\xi$ is obtained as follows. Recall that we were originally dealing with 
bundles $E_0$ with $Sl(2,\bbc)$-structure, and so a map $\xi_0:H^0( \Sigma,
\Lambda^2(E_0))\rightarrow
H^0(\Sigma, {\cal O})=\bbc$. The bundle $E$ is obtained as $E_0\otimes L$, giving 
the induced map $\xi$.
We have a commutative diagram
$$\matrix {H^0(\Sigma,\Lambda^2(E_0))&{\buildrel \xi_0\over\rightarrow}&H^0(\Sigma, {\cal O})\cr
\downarrow ev_{\Lambda^2(E_0)}&&\downarrow ev_{\cal O}\cr
\Lambda^2(E_0)|_{p_i}&{\buildrel \xi^{p_i}_0\over\rightarrow}&\bbc}$$ 
One then has that $\xi_0= ev_{\cal O}^{-1}\circ \xi^{p_i}_0 \circ ev_{\Lambda^2(E_0)}$;
and it is this composition that one uses to extend the definition of $\hat \beta $
to the cases when the bundle acquires torsion at the $p_i$. One obtains then 
maps $\hat \beta_i$, and one checks (by a calculation in local coordinates, say for 
a one-parameter family; see lemma (4.10) above) that $\hat \beta_i$ 
vanishes when the bundle acquires torsion at the $p_i$; when there is no torsion,
 one simply 
has the preceding definition. One then sets 
$$\hat E_X' =( (\hat \beta_1 ,\hat \alpha_1^2),...,(\hat \beta_n,
\hat \alpha_n^2))$$

We  analyse the stability of elements of $Y$.
We had that the presence of torsion in the kernel of $\alpha_i$ destabilised 
an element of $X$;  the same is true of elements of $Y$:

{\sc Lemma (4.12)}: \tensl A semi-stable element $y$ of $Y$ corresponds to a bundle
 $E$ without 
torsion in the kernel of $\alpha$. \tenrm

{\sc Proof}: Let $U = H^0({\cal O}^N)$, and split $U$ as $U_0\oplus U_1$, with $U_0$
the sections of the torsion subsheaf of $E$, and $U_1$ the sections of $E/$torsion. 
The element of $Z_1$ corresponding to  $E$ is represented by an element of 
$V_1$; if in addition there is a non-trivial subspace $W$
of $U_0$ lying in the kernel of $\alpha$, one can destabilise $y$ by a one parameter 
multiplicative subgroup acting with positive weight on $W$ and negative weight 
on the complement.

It follows then that the semi-stable $y$ have their torsion at the $p_i$,
 and that 
it is either $0$ or $\bbc_{p_i}$. From now on we restrict ourselves to sheaves 
whose torsion is of this type.

Now we consider  pairs $((E,\alpha)\in X, \phi)$. As noted above, the sheaf $E$ correponds to an element
of $Z_1$, and the forms $\alpha_i$ gives us non-zero elements
 $\tilde\alpha_i$ of $V_2^*$.
 One can then map our pairs into  
$Y$ by $(E,\alpha, \phi)\mapsto ((\phi(p_1)^N, \phi(p_1)^{N-1}\tilde\alpha_1^2),$
 $(\phi(p_2)^N, \phi(p_2)^{N-1}
\tilde\alpha_2^2),...,$ $(\phi(p_n)^N, \phi(p_n)^{N-1}\tilde\alpha_N^2)$. Now we take the closure $\hat Y$ of 
the image  of this map. This closure allows
the $\alpha_i$ to go to zero (while $\phi(p_i)\ne 0$), while preserving the information of a subspace of $E$ at $p_i$.
Note that $E$ is then torsion free at $p_i$, since $\phi(p_i)\ne 0$.

{\sc Proposition (4.13)}: \tensl Let 
$$y = ((a_{0,1}, a_{1,1}),(a_{0,2}, a_{1,2}),...,(a_{0,n}, a_{1,n}))$$
be an element of $\hat Y$. Let $\Gamma(y)$ be the set of $\gamma_i$ such that $\gamma_i\in [0,1]$, $\gamma_i=0$
if $a_{1,i}=0$, and $\gamma_i=1$ if $a_{0,i}=0$. Then $y$ is semi-stable iff
for one element $\gamma$ of $\Gamma(y)$, $\pi(y)\in X$ is $\gamma$-semi-stable.\tenrm

{\sc Proof}: If $\pi(y)$ is (semi)-stable for one choice of weights in $\Gamma(y)$,
then $y$ is (semi)-stable. This follows from the fact that the $Sl(N,\bbc)$ quotient
of $X$ is just the $SL(N,\bbc)\times (\bbc^*)^n$ quotient of $Y$, with the $\gamma_i$
defining the weights of the $(\bbc^*)^n$ action (Thaddeus[Th]).

On the other hand, let $\pi(y)$ be unstable  for all $\gamma_i$ in $\Gamma(y)$.
We first note that for $a_{1,i}=0$ (hence $\gamma_i=0$), one can simply remove the point 
$p_i$, as it does not affect the definition of stability either in $X$ or in $\hat Y$.
Similarily, if $a_{0,i}=0$, (hence $\gamma_i=1$) one has torsion. Stability for $E$ with 
$\gamma_i=1$ is then equivalent to stability for $E/$torsion with $\gamma_i=0$, and by the preceding
argument, we can drop the point, at the price of considering $E/$torsion instead.

We are then in the situation of having a bundle $E$ of degree $k<0$ which can be destabilised 
by any choice of $\gamma_i\in [0,1]$. Let $\gamma_i = 0$ to start. There is then a subline bundle
of maximal  degree $j$, with $2j>k$. If $L\subset ker(\alpha_i)$, increasing $\gamma_i$ only makes $L$
more unstable, so we drop those points, and assume that $L$ does not lie in any of the kernels.
Now $L$ will destabilise for a choice of $\gamma_i$ if $2j - \Sigma_i\gamma_i > k$; assume that one 
can increase the $\gamma_i$ so that $2j - \Sigma_i\gamma_i = k$, so that $L$ no longer destabilises.
By hypothesis there is another $L'$ which now destabilises. $L'$ cannot be a subsheaf of 
$L$, since then $L$ would still destabilise. This tells us that the degree $j'$ of $L'$ is less
than the degree of $E/L$, that is, $k-j$. This gives
$$k< 2j' +\sum_i(\pm\gamma_i)\leq 2k-2j +\sum_i\gamma_i =2k,\eqno (4.14)$$
a contradiction since $k<0$.

The same line bundle $L$ then destabilises uniformly in the $\gamma_i$,
and so one can construct a one parameter subgroup destabilising 
the corresponding element of $\hat Y$.

Semi-stable elements of $\hat Y$ thus all live above semi-stable elements of $X$
for some choice of weights, and  the quotient $\hat Y//SL(N,\bbc)$ corresponds to 
(S-equivalence classes of) quadruples $(E,\alpha_i,\hat \alpha_i, \phi)$, where $E$ is a rank 2 sheaf, 
$\alpha$ is the map $E\rightarrow \oplus_i\bbc_{p_i}$, $\hat \alpha_i$ is a subspace
of $E|_{p_i}$ which is the kernel of $\alpha$ when $\alpha$ is non-zero, one-dimensional 
when $E$ has no torsion at $p_i$; 
$\phi$ is an $Sl(2,\bbc)$-structure. The pair $(E,\alpha)$ will be semi-stable for
 some choice of weights satisfying  the constraint that $\gamma_i = 0$ when $\alpha_i=0$,
and $\gamma_i=1$ when $E$ has torsion at $p_i$. 

This is not quite the moduli space that we want, as  the 
$\hat \alpha_i$ constitute extra information, but  only when $\alpha_i=0$. At these
points, the $\hat \alpha_i$ lie in $\bbp_1=\bbp(E|_{p_i})$.
Embedding $V_1$ into $W_{1,i}= V_1^{\otimes p}$, a non-zero  element of $L_{0,i}$ can be thought of 
as an element of $W_{1,i}^*$; similarly, a non-zero element of $L_{1,i}$ corresponds to 
an element of $W_{2,i} = V_1^{\otimes p-1}\otimes V_2^{\otimes 2}$; the subvariety $\hat Y$
maps to a subvariety $\tilde Y$ in $\bbp(W_{1,1}\oplus W_{2,1})\times \bbp(W_{1,2}\oplus W_{2,2})
\times...\times\bbp(W_{1,n}\oplus W_{2,n})$. This map collapses the unwanted $\bbp_1$'s, and
otherwise is an embedding. The stability analysis goes through unchanged, since 
for $\gamma_i=0$, the group action does not ``see'' the extra $\bbp_1$, which lives in the 
$Z_2$ factor. 

Let ${\cal P} = \tilde Y//Sl(N,\bbc)$ be the geometric quotient; we will 
call it the moduli space of framed parabolic bundles.

{\bf 4.2 The map from $P$ to ${\cal P}$}

Here again, we only give the map for the case $SU(2)$.

Let $\Sigma$ be a closed Riemann surface, and let $\Sigma_0=
\Sigma-\{p_1,..,p_n\}$. We define a map from the space $P(\Sigma_0)$
to the space of framed parabolic sheaves over $\Sigma$ with 
parabolic structure at the $p_i$.

We have  fixed a rank two ${\cal C}^\infty$ 
bundle
$E$ over $\Sigma$ with $c_1(E) = 0$.
Let $\rho$ be an element of $P$. First of all, it defines a 
representation of the fundamental group of $\Sigma_0$ into 
$SU(2)$, and so a local system over $\Sigma_0$, which then defines 
a holomorphic structure over $\Sigma_0$. Recall that we had points
$\tilde p_i$ near the $p_i$, and choose parametrisations  of 
neighbourhoods $D_i$ of these points in 
$\Sigma$ so that $p_i$ corresponds to $z=0$ and $\tilde p_i$ to $z=1$.
We now extend the holomorphic structure to the $D_i$ as follows.
One has that the monodromy of the local system at the i-th puncture
is given by 
$$\pmatrix {e^{\pi \gamma_i}&0\cr 0&e^{-\pi \gamma_i}}, \gamma_i\in [0,1]\eqno (4.15)$$
Let us first consider the case of $\gamma_i\in (0,1)$. One then has associated 
to each $p_i$ a volume form on the eigenspace, or, what is equivalent,
a trivialisation $t$ in which the monodromy of the puncture is diagonal.
Choose the flat trivialisation on the complement of $z=0$ which 
is equal to $t$ at $\tilde p_i$, we  glue our bundle 
over the punctured disk $D_i^*$ to the trivial bundle 
${\cal O}\oplus {\cal O}$ over $D_i$ by the transition matrix:
$$\pmatrix {z^{ \gamma_i/2}&0\cr 0&z^{ -\gamma_i/2}}.\eqno (4.16)$$
We take as extra structure on the resulting bundle $E$ at $p_i$
the map $v_i:E_{p_i}\rightarrow \bbc$ given in these trivialisations by 
$v_i= ( sin(\pi\gamma_i/2),0)$.

When $\gamma_i$ moves to zero, the extra data in $P$ which gave
 the trivialisation  at $p_i$ disappears; on the other hand, 
the glueing matrix becomes the identity matrix, and $v_i$ becomes zero,
 so the framing is no longer required, and 
the framed parabolic sheaf one obtains is simply the trivial one. 

The case when $\gamma_i$ moves to 1 is more interesting. What we glue in 
over $D_i$ is not a bundle, but the sheaf ${\cal O}\oplus {\cal O}\oplus 
\bbc_0$ which has some torsion at the origin. Away from the origin this is just 
${\cal O}\oplus {\cal O}$, and we have the map
$${\cal O}\oplus {\cal O} {\buildrel \pmatrix {z^{-1}&0\cr0&1}\over \longrightarrow}
 {\cal O}\oplus {\cal O}.\eqno(4.17) $$
If on the left hand side, we have over the disk a volume form $1$
and a map at the origin from the fiber to $\bbc$ given by $(0,1)$, 
we have correspondingly on 
 the right hand side the volume form $z$, and a 
form $V= (0,0) $.
When one composes this with the glueing matrix (4.16) for $\gamma_i=1$,
 one obtains:
$$\pmatrix {z^{ -1/2}&0\cr 0&z^{ -1/2}}.\eqno (4.18)$$
which is central for all values of $z$. Choosing a trivialisation at 
$\tilde p_i$, we use (4.18) to 
glue in the sheaf ${\cal O}\oplus {\cal O}\oplus \bbc_0$, with 
its degenerating $Sl(2,\bbc)$ structure to  our bundle over the
punctured surface, and a $v$ at $z=0$ given by $(0,0,1)$. Because (4.18) is central, the result is independent of the framing
chosen, and so, even though the framing in $P$ disappears at $\gamma_i = 1$,
there is still a well defined map $\rho: P\rightarrow {\cal P}$

{\sc Theorem (4.19)}: \tensl The map $P\rightarrow {\cal P}$ is an isomorphism.

\tenrm {\sc Proof}: We have the natural map:
$$\Gamma: P\rightarrow [0,1]^n$$
which associates to each representation the $(\gamma_1,...,\gamma_n)$
of holonomies at the punctures.
For each value of the $\gamma_i$, we have the result, 
due to Mehta and Seshadri [MS] and Narasimhan and Seshadri [NS], that the moduli 
space of representations with holonomy given by the $\delta_i$ is equivalent
to the moduli space of parabolic bundles with polarisations again determined
by the $\gamma_i$. This tells us, in essence, that ``fiberwise'', our theorem is
already proven.

We construct a map $\Delta:{\cal P}\rightarrow [0,1]^n$
with $\Delta \circ \rho = \Gamma$. Recall that ${\cal P}$ is obtained 
by taking a geometric $SL(N,\bbc)$ quotient of the $(\bbp^1)^n$-bundle 
over $X$: $Y= \bbp(L_{0,1}\oplus L_{1,1})\oplus \bbp(L_{0,2}\oplus L_{1,2})
\oplus\cdots
\bbp(L_{0,n}\oplus L_{1,n})$, and then collapsing some 
$\bbp^1$'s. This collapsing is essentially irrelevant to the present proof,
and so we omit it from now on.  We now recall from [Th] that
 the moduli spaces for 
each set of $\delta_i$ is obtained as a $(\bbc^*)^n\times SL(N,\bbc)$-quotient
of $Y$, where the $\delta_i$ are weights of the 
$(\bbc^*)^n$-action. Now recall (see e.g. [MFK]) the equivalence due to 
Mumford, Guillemin-Sternberg, etc. between symplectic quotients 
and geometric invariant theory quotients. We have, for the Kahler structure under the natural
embedding of $Y$ into projective space given by the ${\cal O}(1,1,...,1))$
polarisation, that the action of $(S^1)^n\subset (\bbc^*)^n$ has a moment map
$$\hat\Delta:Y\rightarrow [0,1]^n.\eqno(4.20)$$
This moment map restricts on each fiber $(\bbp^1)^n$ to the 
standard moment map. Since the action of $(S^1)^n$ commutes with that of $SL(N,\bbc)$, this map
is  $SL(N,\bbc)$-invariant, and so  $\hat \Delta$ descends to $P$,
and gives a map $\Delta$. To see that it has the right properties, note that
$$\eqalign{Y//(\bbc^*)^n\times SL(N,\bbc) &\simeq 
[\hat\Delta^{-1}(\delta_1,...,\delta_n)/(S^1)^n]//SL(N,\bbc)\cr
&\simeq \Delta^{-1}(\delta_1,...,\delta_n)/(S^1)^n}\eqno (4.21)$$
for each parabolic moduli space. From this, one obtains  that 
$\Gamma =\Delta\circ \rho$. As the map from $P$ to ${\cal P}$ commutes with the 
$(S^1)^n$-action, and we have isomorphisms $\Gamma^{-1}(\gamma)/(S^1)^n\simeq
 \Delta^{-1}(\gamma)/(S^1)^n$, we are done.

{\bf 4.3 The case of a trinion.}

We now consider the special case of $\bbp^1$ with three marked points.
We know from the results above that the moduli space is $\bbp^3(\bbc)$,
but it is instructive to obtain this as a complex quotient, though the quotient
here is not the standard one.  We first
examine which sheaves $E$ on $\bbp_1$ arise from semi-stable framed 
parabolic bundles:

When $E$ is  a bundle,  the semi-stability condition tells that either 
$E\simeq {\cal O}\oplus {\cal O}$ or ${\cal O}(1)\oplus {\cal O}(-1)$,
for if $E\simeq {\cal O}(j)\oplus {\cal O}(-j), j\ge 2$, one finds that 
${\cal O}(j)$ is a destabilising subbundle, no matter what the map 
$\alpha$ is. When the bundle is trivial, one finds that the framed bundle
 is stable
when none of the $\alpha_i$ are non-zero and
 any two of the subspaces $ker(\alpha_i), i=1,2,3$ span $\bbc^2$ in the global trivialisation.

When $E$ has torsion at one point $p$, the bundle $E^{**}=E$/torsion must be
of the form ${\cal O}(-1)\oplus {\cal O}$. The other possibilities would 
be ${\cal O}(-j)\oplus {\cal O}(j-1), j\ge 2$, and then the sheaf
${\cal O}(j-1) \oplus \bbc_p$ would be destabilising.

Similarily, when $E$ has torsion at two points, one finds that $E\simeq 
{\cal O}(-1)\oplus {\cal O}(-1)\oplus \bbc_p\oplus\bbc_q$. For $E$
to be semistable, $E$ cannot have three torsion points. 

This discussion shows that the sheaves corresponding to semi-stable points 
are all such that $E(1)$ is generated by four global sections, 
from which it follows:

{\sc Proposition (4.22)}: \tensl  $E$ fits into an exact sequence:
$${\cal O}(-2)^{\oplus 2}{\buildrel {A+Bz}\over\longrightarrow}
{\cal O}(-1)^{\oplus 4}\rightarrow E\rightarrow 0,\eqno (4.23)$$
 and the map ${\cal O}(-2)^{\oplus 2}\longrightarrow
{\cal O}(-1)^{\oplus 4}$ is injective (as a map of bundles) 
away from the torsion points of $E$.\tenrm

The bundle is then determined by two $4\times 2$ matrices.
 The $v_i, i=1,2,3$ can then be realised as maps ${\cal O}^{\oplus 4}(-1)|_{p_i}
\rightarrow \bbc$ which annihilate the image of ${\cal O}(-2)^{\oplus 2}$; 
letting the points $p_i$ correspond to $z=0,1,\infty$, this gives us 
row vectors $V_1$ in $\bbc^4$ satisfying 
the conditions
$$V_1A=0,\  V_2(A+B)=0, \  V_3B=0\eqno (4.24)$$
The isomorphism $\Lambda^2(E) = (\Lambda^2{\cal O}(-2)^{\oplus 2})^*
\otimes
\Lambda^4({\cal O}(-1)^{\oplus 4})$ defines a volume form on $E$, and so
we have a framed parabolic $SL(2,\bbc)$-bundle.

The group $S(Gl(2,\bbc)\times Gl(4,\bbc))$ acts on this data, preserving the 
isomorphism class, by:
$$(g,G)(A,B,V_i)= (GAg^{-1}, GBg^{-1}, V_iG^{-1}),\eqno (4.25)$$
and so we are really interested in the orbit space of this action.
There is a four-tuple of functions which are projectively invariant 
under the action:
$$det (A,B), det \pmatrix{V_1B \cr V_2A}, det \pmatrix{V_1B \cr V_3A}, 
det \pmatrix{V_2A \cr V_3B}.\eqno (4.256)$$
We also note that if $D$ is a multiple of the $2\times 2$ identity matrix, 
then the action of $(D^2, diag(D,D))$ simply rescales the data 
$(A,B, V_i)$. We can then simply consider the orbits in the projectivisation
of the  space of $(A,B,V_i)$ under the action of the smaller group 
$Sl(2,\bbc)\times Sl(4,\bbc).$

{\sc Proposition (4.27)}:
\tensl The moduli space of semi-stable framed parabolic sheaves
is the geometric quotient of the subset $V$ of 
$\bbp(M(2,4)^2\times (\bbc^4)^3)$ cut out by conditions (4.24), by the group
$Sl(2,\bbc)\times Sl(4,\bbc)$. \tenrm

{\sc Lemma (4.28)}: \tensl A point $(A,B,V_i)$ is semi-stable for this action if and only if
one of the invariants of (2) is non-zero. \tenrm

{\sc Proof}: For an element whose invariants do not vanish, the 
orbit is necessarily semi-stable. Conversely, for an orbit whose invariants 
are zero it suffices to find a destabilising $\bbc^*$ in $Sl(2,\bbc)\times
Sl(4,\bbc)$, that is an whose action on $(A,B,V_i)$ gives 
an orbit with zero in its closure. The fact that the invariants all vanish 
tells us that we can find a basis in $\bbc^2$ such that 
$V_iA,V_iB$ are all of the form $(*,0)$. 
We subdivide into cases, according to the 
size of the span in $\bbc^4$ of the $V_i$. For the case when the span is 
three-dimensional, one sets $V_i= (*,*,*,0)$, and one then has that both $A$ and $B$ are of the form 
$$\pmatrix{*&0\cr *&0\cr *&0\cr  *&*}.\eqno (4.29)$$
Now choosing integers $n,m$ such that $0<n<m<3n$, we find that 
the subgroup of elements $(diag (z^m, z^{-m}), diag(z^{-n},z^{-n},z^{-n},z^{3n})),
z\in \bbc$ destabilises. The other cases are similar, if slightly more elaborate.

We then have a well-defined map $\Phi$ from the moduli of semi-stable framed parabolics to 
$\bbp^3(\bbc)$. 

{\sc Theorem (4.30)} \tensl The map $\Phi$ is an isomorphism.\tenrm

{\sc Proof}: We must check that the closed orbits map bijectively into 
$\bbp^3(\bbc)$. Again, there is a case by case analysis.

{\it Case (i)} $Det(A,B)\ne 0$. One can choose bases so that 
$$A= \pmatrix {1&0\cr 0&1\cr 0&0\cr 0&0}\ B=\pmatrix{0&0\cr 0&0\cr1&0\cr 0&1}\eqno (4.31)$$
The stabiliser of this form is $Sl(2,\bbc)$. Conditions (4.24) tell us that
$$ V_1= (0,0,a,b), V_2 = (c,d, -c,-d), V_3 = (e,f,0,0).\eqno (4.32)$$
One is left with the problem of quotienting out the action of $SL(2,\bbc)$
on the three-tuple $(a,b), (c,d), (e,f)$  of vectors in $\bbc^2$. The stable
 points with finite stabiliser are those for which the three vectors span $\bbc^2$.
The three functions 
$$det \pmatrix{V_1B \cr V_2A}, det \pmatrix{V_1B \cr V_3A}, 
det \pmatrix{V_2A \cr V_3B}\eqno (4.33)$$ 
map this set bijectively to $\bbc^3\backslash\{0\}$.
All the others are just semi-stable, with the single closed 
orbit $((0,0),(0,0),(0,0))$, which corresponds to the origin in $\bbc^3$. 
All of these orbits have coordinate $(1,0,0,0)$ in $\bbp^3$.

{\it Case (ii)} $Det(A,B) = 0 , Im(A)+Im(B)$ three dimensional, 
and $(Az_0+Bz_1)$ injective for all $(z_0,z_1)\ne (0,0)$: One can then normalise
$(A,B)$ to 
$$A= \pmatrix {1&0\cr 0&1\cr 0&0\cr 0&0}\ B=\pmatrix{0&0\cr 1&0\cr0&1\cr 0&0}\eqno (4.34)$$
One then finds that 
$$ V_1= (0,0,a,b), V_2 = (c,-c, c,d), V_3 = (e,0,0,f).\eqno (4.35)$$
The projective coordinate functions of (4.26) are $(0,-ac, -ae, ec)$.
When these are non-zero, the stabiliser is finite and the points are 
stable. When one of $a,c,e$ is zero, then the orbit is not closed, and its 
coordinates correspond to one of the points $(0,1,0,0), (0,0,1,0)$ or $(0,0,0,1)$
in $\bbp^3$.

{\it Case (iii)} $Det(A,B) = 0 , Im(A)+Im(B)$ three dimensional, 
but  $(Az_0+Bz_1)$ fails to be injective for some $(z_0,z_1)\ne (0,0)$:
One checks that the stability condition forces $z_0/z_1$ to be either
$0,1$ or $\infty$, and that injectivity fails at only one point.
 We will suppose that the point is zero; the other cases are similar.
We can then normalise $A,B$ to:
$$A= \pmatrix {0&1\cr 0&0\cr 0&0\cr 0&0}\ B=\pmatrix{0&0\cr 1&0\cr0&1\cr 0&0}.\eqno (4.36)$$
One then has that 
$$V_1=(0,a,b,g), V_2=(c,0,-c,d), V_3= (e,0,0,f)\eqno (4.37)$$
and the projective coordinate functions of (1) are
$(0,ac, ae, 0)$. From this, one has that $a$ must be non-zero, and that 
either $c$ or $e$ is non-vanishing. When both of them are non-zero, the 
point is stable, and when they are not, the orbit is not closed.

{\it Case (iv)} $Det(A,B) = 0 , Im(A)+Im(B)$ two dimensional: Here one finds that
$(Az_0+Bz_1)$ fails to be injective for 
two elements $z_0/z_1$ of $\bbp^1$, and that these points must lie in $0,1,\infty$.
Again, we do one case, that of the points $0,\infty$, the others being similar.
One can normalise to 
$$A= \pmatrix {0&1\cr 0&0\cr 0&0\cr 0&0}\ B=\pmatrix{0&0\cr 1&0\cr0&0\cr 0&0}.\eqno (4.38)$$
 with 
$$V_1=(0,a,b,c), V_2=(0,0,d,e), V_3= (f,0,g,h)\eqno (4.39)$$
and projective  coordinates $(0,0,-ag,0)$. The unique closed orbit is that 
of 
$$V_1=(0,a,0,0), V_2=(0,0,0,0), V_3= (f,0,0,0)\eqno (4.40)$$

This in fact exhausts all cases of semi-stable orbits. Recapitulating, we see 
that all the points of $\bbp^3$ apart from $(1,0,0,0), (0,1,0,0), (0,0,1,0)$ 
and $(0,0,0,1)$ correspond to stable points, with the four exceptional 
points each corresponding to a unique closed semi-stable orbit. 

The four cases given above each correspond to different types of sheaf.
In case (i), the condition $Det(A,B)\ne 0$ corresponds to the sheaf being
${\cal O}\oplus {\cal O}$, and one has a $\bbc^3$ of framed parabolics with 
such a structure. In case (ii), one has a bundle ${\cal O}(-1)\oplus
{\cal O}(1)$, and the space of framed parabolics with such a structure is 
the plane $\bbp^2$ at infinity, minus the three coordinate lines.
In case (iii), there is one point $p$  amongst $0,1,\infty$ at which there is torsion,
and the sheaf is of the form ${\cal O}\oplus {\cal O}(-1)\oplus \bbc_p$. The
parabolics corresponding to these correspond to the three coordinate lines
in the plane at infinity, minus their intersections. Finally, the points
in case (iv) correspond to sheaves ${\cal O}(-1)\oplus {\cal O}(-1)\oplus 
\bbc_q\oplus\bbc_p$, with $\{p,q\}\in \{0,1,\infty\}$: there are three choices
and a corresponding three points in the plane at infinity, given by the 
intersection of the coordinate axes.

{\bf 4.4. A glueing map and a toric variety}

Given a pair of Riemann surfaces $\Sigma_0,\Sigma_1$ with marked
points $p_0,p_1$, one can define a singular nodal curve $\Sigma$
by identifying the points $p_0$ and $p_1$. There is also on the 
level of the moduli spaces ${\cal P}_0, {\cal P}_1$ of framed
parabolic sheaves a (partial) glueing map which one can define as follows:
one has for a framed parabolic sheaf $E_0$ over $\Sigma_0$ 
a map $\alpha_0:(E_0)_{p_0}\rightarrow \bbc$,
and similarly for $\Sigma_1$. One can combine the two into a diagram:
$$ (E_0)_{p_0}\rightarrow \bbc\leftarrow (E_1)_{p_1}.$$
If we consider two such diagrams to be equivalent if   the framed parabolic bundles
on both $\Sigma_0$ and $\Sigma_1$ are isomorphic and if the induced maps
$(E_0)_{p_0}/ker(\alpha_0)\rightarrow (E_1)_{p_1}/ker(\alpha_1)$ are the same,
then taking equivalence classes of such diagrams amounts to quotienting 
out the anti-diagonal action of $\bbc^*$ on the framings at $p_0$, $p_1$,
and we define the moduli space of framed parabolic sheaves over 
$\Sigma$ to be the geometric quotient:
$${\cal P}_{\Sigma}={\cal P}_0\times  {\cal P}_1//\bbc^*.$$
This is the complex analogue of the symplectic glueing defined above.
>From the symplectic point of view, the space $P$ which was very close to the 
moduli space $\mT_{\tilde \Sigma}$ of the glueing of $\Sigma_0$ and $\Sigma_1$ into a smooth 
surface $\tilde \Sigma$.
In this picture, the geometry of what the glued moduli corresponds to 
is more evident: one is dealing with a moduli space of a degenerate curve,
and one should be able to obtain ${\cal P}_{\Sigma}$ as a degeneration of $\mT_{\tilde \Sigma}$.

We note that there still remains the diagonal action of $\bbc^*$ at
the puncture.  Now let us consider a decomposition of
a smooth surface $\tilde \Sigma$ into trinions, and let us pinch the boundaries of the 
trinions to points, so that one has a singular surface $\Sigma$, consisting of
$2g-2$ spheres touching at $3g-3$ points.  The moduli space ${\cal P}_{\Sigma}$ 
is $3g-3$ dimensional and is obtained 
from the glueing of $2g-2$ trinion moduli spaces ``at'' $3g-3$ points. It has an 
action of $\bbc^*)^{3g-3}$, and is a toric variety.

For $\bbp^3$, the moduli space of the trinion, 
the action of $(a,b,c)\in(C^*)^3$ is  
 $(x,y,z,w)\mapsto (x, aby, ac z, bcw)$ in the coordinates given 
in the previous section. Its moment map is given by the holonomies, and 
the moment polytope is the tetrahedron of section 3.
 The moment polytope $PT_\Sigma$
of ${\cal P}_{\Sigma}$ is given by taking the intersection with some hyperplanes of the
product of $2g-2$ tetrahedra.

The main purpose of [JW] is to provide a   justification 
for the Verlinde formulae in terms of real quantisation. To do this, 
as we saw, they exhibited a toric structure associated to a trinion decomposition of the
surface, which was however ill defined when the holonomy of the flat connections
along the circles associated to the decomposition was central. There is still, however,
a well defined moment map to a polytope. This polytope, as one can check in [JW],
is exactly $PT_{\Sigma}$.

Referring to the results of section 3, one therefore has a commuting diagram
$$\matrix{\mT_{\tilde \Sigma}&&\rightarrow&&{\cal P}_{\Sigma}\cr
&\searrow&&\swarrow&\cr
&&PT_{\Sigma}&&}$$

On one hand, for $M_{\tilde \Sigma}$, the dimension of the spaces of sections of 
line bundles over $M_{\tilde \Sigma}$ (which is what the Verlinde formula computes)
is given by the heuristic of Bohr-Sommerfeld quantisation in terms of the number
of (fractional) integer points in the polytope $PT_{\Sigma}$. On the other,
for ${\cal P}_{\Sigma}$,
the theory of toric varieties tells us that the space of sections 
of the corresponding line bundle over ${\cal P}_{\Sigma}$, is given in terms of the same 
count of points. Furthermore, the correspondence between $\tilde \Sigma$ and $\Sigma$  
is that of a degeneration of the surface; this leads us to believe that a symplectic 
cobordism can be built between $M_{\tilde \Sigma}$ and ${\cal P}_{\Sigma}$, which 
preserves the dimension of the space of sections, giving (yet another) 
proof of the Verlinde formula.

\bigskip
\parindent = 0pt
{\bf Bibliography}
\bk
 
[AMM] A. Alekseev, A. Malkin, E. Meinrenken, Lie group valued
moment maps. {\it J. Diff. Geom.} {\bf 48} (1998) 445-495.

[AB] M.F. Atiyah and R. Bott, The Yang-Mills equations over 
Riemann surfaces,  {\it Phil. Trans. Roy. Soc. Lond.} {\bf A308},
(1982) 523-615. 

[Bh] U. Bhosle, Parabolic vector bundles on curves, {\it Ark. Mat.}
{\bf 27}, 15-22 (1989).

[D] T. Delzant, Hamiltoniens p\'eriodiques et images convexes 
de l'application moment, {\it Bull. Soc. Math. France} {\bf 116}
(1988) 315-339.

[Gi] D. Gieseker, On the moduli of vector bundles on algebraic
surfaces, {\it Ann. Math.} {\bf 106}, 45-60 (1977).

[G1] W. Goldman, Invariant functions on Lie groups and Hamiltonian
flows of surface group representations, {\it Invent. Math.} 
{\bf 85}, 263-302 (1986).

[G2] W. Goldman, The symplectic nature of fundamental groups of 
surfaces.  {\it Advances in Math.} {\bf 54}, 200-225 (1980).

[GS] V. Guillemin and S. Sternberg, {\it Symplectic Techniques in Physics},
Cambridge University Press, 1984.

[H] J. Huebschmann, On the variation of the Poisson structures of 
certain moduli spaces, preprint dg-ga/9710033.

[HL] D. Huybrechts, M. Lehn, Stable pairs on curves and surfaces.
{\it J. Algebraic Geom.} {\bf 4}, 67-104 (1995).

[J1] L.C. Jeffrey, Extended moduli spaces of flat connections 
on Riemann surfaces. {\it Math. Annalen} {\bf 298},   667-692 (1994).
%corresponds to one reference to [Jeffrey] and several to [J]

[J2] L.C. Jeffrey, Symplectic forms on moduli spaces of flat
connections on 2-manifolds. 
In  Proceedings of the Georgia
International Topology Conference (Athens, GA, 1993), ed. W. Kazez,
Amer. Math. Soc./International Press
 AMS/IP Studies in Advanced Mathematics
{\bf 2}, 268-281 (1997).

[JW] L.C. Jeffrey and J. Weitsman, Bohr-Sommerfeld orbits in 
the moduli space of flat connections and the Verlinde dimension
formula. {\it Commun. Math. Phys.} {\bf 150}, 593-630 (1992). 

[MFK] D. Mumford, J. Fogarty, F. Kirwan, {\it Geometric 
Invariant Theory}, Springer-Verlag, 1994, chap. 8.2.

[MS] V. Mehta, C.S. Seshadri, Moduli of vector bundles on 
curves with parabolic structure, {\it Math. Ann.}
{\bf 248}, 205-239 (1980).

[MW] E. Meinrenken, C. Woodward, A symplectic proof of Verlinde
factorization, {\it J. Diff. Geom.}, to appear.

[NS] M.S. Narasimhan, C.S. Seshadri, Stable and unitary vector bundles
on a compact Riemann surface.
{\it Annals of Math.} {\bf 82}, 540-567.

[Sj] R. Sjamaar, Imploded cross-sections, in preparation.

[Th] M. Thaddeus, Geometric invariant theory and flips, 
{\it J. Amer. Math. Soc.} {\bf 9}, 691-723 (1996).

\ninerm
\bigskip
\bigskip

\line{   J. C. Hurtubise \hfil L.C. Jeffrey}
\line{Department of Mathematics and Statistics \hfil
Department of Mathematics}
\line{ McGill University  \hfil University of Toronto}
\line{email: hurtubis@math.mcgill.ca \hfil email: jeffrey@math.utoronto.ca}

\end